\documentclass[twoside,leqno]{article}
\usepackage{ltexpprt}
\usepackage{url}
\usepackage{stackengine}

\usepackage{booktabs}
\usepackage{amsmath}
\usepackage{amsfonts}
\usepackage{mathtools}
\usepackage{cleveref}
\usepackage{xcolor}
\usepackage{subcaption}
\usepackage{algorithm,algpseudocode}
\crefname{equation}{Eq.}{Eqs.}
\crefname{figure}{Fig.}{Figs.}
\newcommand{\grad}{\nabla}

\newcommand{\comment}[1]{}
\DeclarePairedDelimiter{\norm}{\lVert}{\rVert}

\AtBeginDocument{
	\abovedisplayskip=1ex
	\belowdisplayskip=1ex
	\abovedisplayshortskip=0ex
	\belowdisplayshortskip=1ex
}

\newcommand{\rev}[1]{#1}

\begin{document}
%%%%%%%%%%%%%%%%%%%%%%%%%%%%%%%%%%%%%%%%%%%%%%%%%
%%%%%%%%%%%%%%%%%%%%%%%%%%%%%%%%%%%%%%%%%%%%%%%%%
%%%%%%%%%%%%%%%%%%%%%%%%%%%%%%%%%%%%%%%%%%%%%%%%%

\title{Nesterov Acceleration of Alternating Least Squares for Canonical Tensor Decomposition: Momentum Step Size Selection and Restart Mechanisms\protect\thanks{This is a revised and extended version of paper arXiv:1810.05846v1 (2018).}}
%\title{Nesterov Acceleration of Alternating Least Squares for Canonical Tensor Decomposition: Step Size Selection and Restart Mechanism\protect\thanks{This is a revised and extended version of paper arXiv:1810.05846 (2018).}}

\author{
Drew Mitchell\thanks{School of Mathematical Sciences, Monash University, demit3@student.monash.edu}
\and
Nan Ye\thanks{School of Mathematics and Physics, University of Queensland, nan.ye@uq.edu.au}
\and
Hans De Sterck\thanks{Department of Applied Mathematics, University of Waterloo, hdesterck@uwaterloo.ca} 
}
\date{}

\maketitle

%\presentaddress{}

\begin{abstract} \small\baselineskip=9pt 
We present Nesterov-type acceleration techniques for Alternating Least Squares
(ALS) methods applied to canonical tensor decomposition.
While Nesterov acceleration turns gradient descent into an optimal first-order
method for convex problems by adding a momentum term with a specific weight
sequence, a direct application of this method and weight sequence to ALS results in erratic
convergence behaviour. This is so because ALS is accelerated instead of gradient descent for our non-convex problem.
Instead, we consider various restart mechanisms and suitable choices of momentum weights
that enable effective acceleration.
Our extensive empirical results show that the Nesterov-accelerated ALS methods with restart can be dramatically more efficient
than the stand-alone ALS or Nesterov accelerated gradient methods, when problems are ill-conditioned or accurate solutions
are desired. The resulting methods perform competitively with or superior to existing acceleration methods for ALS, including
ALS acceleration by NCG, NGMRES, or LBFGS, and additionally enjoy the benefit of being much easier to implement. We also compare with Nesterov-type updates where the
momentum weight is determined by a line search, which are equivalent or closely related to existing line search methods for ALS.
On a large and ill-conditioned 71$\times$1000$\times$900 tensor consisting of readings from chemical sensors to track hazardous gases, the restarted Nesterov-ALS method \rev{shows desirable robustness properties and} outperforms any of the existing methods we compare with
by a large factor. 
There is clear potential for extending our Nesterov-type acceleration approach to accelerating other optimization algorithms
than ALS applied to other non-convex problems, such as Tucker tensor decomposition.
Our Matlab code is available at
\url{https://github.com/hansdesterck/nonlinear-preconditioning-for-optimization}.
\end{abstract}


%\footnotetext{test}

%%%%%%%%%%%%%%%%%%%%%%%%%%%%%%%%
\section{Introduction}
%%%%%%%%%%%%%%%%%%%%%%%%%%%%%%%%

Nesterov's accelerated gradient descent method is a celebrated method for
speeding up the convergence rate of gradient descent,
achieving the optimal convergence rate obtainable for first order methods
on convex problems \cite{nesterov1983method}. 
Nesterov's method is an extrapolation method that specifies carefully tuned
weight sequences for the so-called momentum term which updates an iterate in
the direction of the previous update. With those tailored weight sequences, optimal
convergence rates can be proved for convex problems.

Recent work has seen extensions of Nesterov's accelerated gradient method in several ways:
either the method is extended to non-convex optimization problems \cite{ghadimi2016accelerated,li2015accelerated},
or Nesterov's approach is applied to accelerate convergence of methods that are not
directly of gradient descent-type, such as the Alternating Direction Method of Multipliers
(ADMM) applied to convex problems\cite{goldstein2014fast}.

%This paper attacks both of these challenges at the same time for the canonical tensor decomposition problem: we develop Nesterov-accelerated algorithms for the non-convex CP tensor decomposition problem, and we do this by accelerating ALS steps instead of gradient descent steps.

In this paper our goal is to extend Nesterov extrapolation to accelerating the Alternating Least Squares (ALS) method for computing the canonical approximation of a tensor by a sum of $R$ rank-one tensors
--- the so-called Canonical Polyadic (CP) decomposition of a tensor\cite{kolda2009tensor}.
As such, this paper attacks two challenges at the same time: we develop Nesterov-accelerated algorithms for a non-convex problem --- the CP tensor decomposition problem ---, and we do this by accelerating ALS steps instead of gradient descent steps.
Since we accelerate ALS instead of gradient descent for our non-convex CP decomposition problem, we cannot simply use the standard sequences for the Nesterov momentum weights that lead to optimal convergence of Nesterov's accelerated gradient method applied to convex problems; in fact, we will illustrate that using these sequences leads to erratic or divergent convergence behaviour. Instead, we investigate in this paper choices for selecting the momentum step size combined with restarting mechanisms that lead to effective acceleration methods for ALS applied to tensor decomposition. Our approach is partially inspired by recent related work on acceleration methods for nonlinear systems\cite{nguyen2018accelerated} and on restarting mechanisms that improve Nesterov convergence for convex problems\cite{odonoghue2015adaptive,su2016differential,goldstein2014fast}.
Our approach and goals are most closely related to (independent) recent work by Ang and Gillis\cite{ang2019accelerating}, who consider Nesterov-type acceleration for Nonnegative Matrix Factorization, and also propose step size and restart mechanisms, which are different from ours.

The Nesterov acceleration methods for ALS with our proposed step size and restart mechanisms are attractive because they are simple to implement, and we show in extensive numerical tests that they are competitive with other recently proposed nonlinear acceleration methods for ALS applied to CP tensor decomposition, most of which are much more involved in terms of implementation. We also believe our step size and restart strategies are of broader interest and may be applied to Nesterov-type acceleration of other simple optimization methods of alternating or (block) coordinate descent-type applied to potentially non-convex problems.

As another contribution of this paper, we also establish links between Nesterov-type acceleration of ALS and other existing acceleration methods for ALS.
We establish links between Nesterov-type extrapolation and the well-known nonlinear acceleration methods of Anderson acceleration and the Nonlinear Generalized Minimal Residual Method (NGMRES) 
\cite{washio1997krylov,oosterlee2000krylov,sterck2012nonlinear,brune2015composing}. In particular, we show that the Nesterov extrapolation formula is equivalent to Anderson acceleration with window size one, and closely related to NGMRES with window size one. In our numerical results, we compare with NGMRES acceleration of ALS.
Furthermore, if the Nesterov momentum weight is determined by a line search, Nesterov acceleration of ALS is closely related to line search methods for ALS that go back to the 1970s and have more recently been enhanced\cite{harshman1970foundations,rajih2008enhanced,chen2011new,sorber2016exact}. We note that acceleration of ALS by NGMRES, NCG and LBFGS also employs line searches\cite{sterck2012nonlinear,sterck2015nonlinearly,sterck2018nonlinearly}. To explore the link with line search methods\cite{harshman1970foundations,rajih2008enhanced,chen2011new,sorber2016exact}, we also compare in our numerical results with Nesterov weights determined by a cubic line search. Line searches are common in acceleration methods for ALS, but are not commonly considered when Nesterov acceleration is applied to optimization algorithms in the literature, so we argue based on our numerical results that Nesterov momentum weights determined by a line search may be a useful algorithmic approach.
Finally, we also explain how Nesterov acceleration of ALS can be interpreted as using ALS as a nonlinear preconditioner for Nesterov's accelerated gradient formula, following a general framework for nonlinear preconditioning of optimization methods that was previously applied to the Nonlinear Conjugate Gradient (NCG) and Limited-Memory Broyden-Fletcher-Goldfarb-Shanno (LBFGS) methods\cite{sterck2018nonlinearly}. 

%\cite{liavas2017nesterov}
%\cite{rajih2008enhanced}
%\cite{sorber2016exact}
%\cite{chen2011new}
%\cite{brune2015composing}
%\cite{fang2009two,walker2011anderson}
%\cite{washio1997krylov,oosterlee2000krylov}
%\cite{ang2019accelerating}

%-----------------------------------------------------------
\subsection{Nesterov's Accelerated Gradient Method.}
\label{subsec:Nesterov}
%-----------------------------------------------------------

Consider the problem of minimizing a function $f(x)$,
\begin{equation}
	\min_{x} f(x).
\end{equation}
Nesterov's accelerated gradient descent starts with an initial guess
$x_{1}$.
For $k \ge 1$, given $x_{k}$, a new iterate $x_{k+1}$ is obtained by first adding a
multiple of the \emph{momentum} $x_{k} - x_{k-1}$ to $x_{k}$ to obtain an
auxiliary variable $y_{k}$, and then performing a gradient descent step at $y_{k}$.
The update equations at iteration $k \ge 1$ are as follows:
\begin{align}
	y_{k} &= x_{k} + \beta_{k} (x_{k} - x_{k-1}),  \label{eq:nesterovy} \\
	x_{k+1} &= y_{k} - \alpha_{k} \grad f(y_{k}), \label{eq:nesterov}
\end{align}
where the gradient descent step length $\alpha_{k}$
and the momentum weight $\beta_{k}$ are suitably chosen numbers, and 
$x_{0} = x_{1}$ so that the first iteration is simply gradient descent.

There are a number of ways to choose the
$\alpha_{k}$ and $\beta_{k}$ so that Nesterov's
accelerated gradient descent converges at the optimal $O(1/k^{2})$ in function
value for smooth convex functions.
For example, when $f(x)$ is a convex function with $L$-Lipschitz gradient, by choosing 
$\alpha_{k} = \frac{1}{L}$, and $\beta_{k}$ as
\begin{align} 
	\lambda_{0} &= 0, \quad
	\lambda_{k} = \frac{1 + \sqrt{1 + 4 \lambda_{k-1}^{2}}}{2}, \label{eq:lambda} \\
	\beta_{k} &= \frac{\lambda_{k-1} - 1}{\lambda_{k}}. \label{eq:momentum1}
\end{align}
one obtains the following $O(1/k^{2})$ convergence rate:
\begin{equation}
	f(x_{k}) - f(x^{*}) \le \frac{2 L \norm{x_1 - x^{*}}}{k^{2}},
\end{equation}
where $x^{*}$ is a minimizer of $f$.
See, e.g., Su et al.\cite{su2016differential} for more discussion on the choices of
momentum weights.
%It adds a \emph{momentum} term to the iterate, and is known to achieve the optimal
%convergence rate achievable by first order methods on convex problems
%\cite{nesterov1983method}.

It will be useful for later considerations in this paper to rewrite the expressions for Nesterov acceleration
in Eqs.\ (\ref{eq:nesterovy}-\ref{eq:nesterov}) solely in terms of $x_k$ or $y_k$.
For the $x_k$ variables, we obtain
\begin{align} 
%        x_{k} &= x_{k-1} + \beta_{k-1} (x_{k-1} - x_{k-2})-\alpha_{k-1} \nabla f \left(x_{k-1} + \beta_{k-1} (x_{k-1} - x_{k-2})\right),
           x_{k} &= q_{\nabla,k-1}\left(y_{k-1}\right),  \nonumber \\
        &=q_{\nabla,k-1}\left(x_{k-1} + \beta_{k-1} (x_{k-1} - x_{k-2})\right),\label{eq:nesterovxonly}
\end{align}
where the function $q_{\nabla,k-1}(\cdot)$ is introduced to represent the result of one steepest-descent step with step size $\alpha_{k-1}$:
\begin{equation}
        q_{\nabla,k-1}(x)=x-\alpha_{k-1} \nabla f(x).  \label{eq:defq}
\end{equation}
Similarly, for the $y_k$ variables, we can write
\begin{align} 
        y_{k} &= q_{\nabla,k-1}(y_{k-1}) + \beta_{k} \left( q_{\nabla,k-1}(y_{k-1}) - q_{\nabla,k-2}(y_{k-2})\right).\label{eq:nesterovyonly}
\end{align}

Recently, there are a number of works that apply Nesterov's acceleration
technique to non-convex problems.
A modified Nesterov accelerated
gradient descent method has been developed that enjoys the same convergence guarantees as
gradient descent on non-convex optimization problems, and maintains the optimal
first order convergence rate on convex problems\cite{ghadimi2016accelerated}.
A Nesterov accelerated proximal gradient algorithm was developed that is guaranteed
to converge to a critical point, and maintains the optimal first order
convergence rate on convex problems\cite{li2015accelerated}.

Nesterov's accelerated gradient method is known to exhibit oscillatory behavior on
convex problems. An interesting discussion on this is provided by Su et al.\cite{su2016differential} 
which formulates an ODE as the continuous time analogue of Nesterov's method.
Such oscillatory behavior happens when the method approaches convergence, and
can be alleviated by restarting the algorithm using the current iterate as the
initial solution, usually resetting the sequence of momentum weights to its
initial state close to 0. An explanation has been provided of why
resetting the momentum weight to a small value is effective using the ODE formulation
of Nesterov's accelerated gradient descent\cite{su2016differential}.
The use of adaptive restarting has been explored
for convex problems,\cite{odonoghue2015adaptive} and Nguyen et al.\cite{nguyen2018accelerated} explored the use of adaptive
restarting and adaptive momentum weights for nonlinear
systems of equations resulting from finite element approximation of PDEs.
Our work is the first study of a general Nesterov-accelerated ALS scheme.

%Several ALS-specific nonlinear acceleration techniques have been developed recently as
%discussed in the introduction \cite{sterck2012nonlinear,sterck2015nonlinearly,sterck2018nonlinearly}.
%These algorithms often have complex forms and incur significant computational
%overhead.
%Our Nesterov-ALS scheme is simple and straightforward to implement, and only
%incurs a small amount of computational overhead.

%-----------------------------------------------------------
\subsection{Nesterov Extrapolation as a Generic Nonlinear Acceleration Scheme} \label{subsec:Nesterov-accel}
%-----------------------------------------------------------
We can now generalize the Nesterov extrapolation approach and apply it to accelerate other optimization methods than steepest descent. Considering a generic iterative optimization method with update formula
$$
x_{k}=q(x_{k-1}),
$$
we replace the steepest-descent operator $q_{\nabla,k-1}(\cdot)$ by the generic update function $q(\cdot)$
in Nesterov update formulas Eq.\ (\ref{eq:nesterovxonly}) or Eq.\ (\ref{eq:nesterovyonly}) to obtain 
\begin{align} 
        x_{k} =q\left(x_{k-1} + \beta_{k-1} (x_{k-1} - x_{k-2})\right),\label{eq:nesterovqxonly}
\end{align}
\begin{align} 
        y_{k} = q(y_{k-1}) + \beta_{k} \left( q(y_{k-1}) - q(y_{k-2})\right).\label{eq:nesterovqyonly}
\end{align}
This extrapolation approach has been used previously to accelerate, for example, ADMM, where $q(\cdot)$ represents an ADMM update\cite{goldstein2014fast}.
Nesterov's technique has also been used to accelerate an approximate Newton method\cite{ye2017nesterov}.

In this paper, we consider Nesterov-type acceleration of ALS for CP tensor decomposition. We write $q_{ALS}(\cdot)$ for the ALS update function, and obtain accelerated formulas
\begin{align} 
        x_{k} =q_{ALS}\left(x_{k-1} + \beta_{k-1} (x_{k-1} - x_{k-2})\right),\label{eq:nesterovALSxonly}
\end{align}
and
\begin{align} 
        y_{k} = q_{ALS}(y_{k-1}) + \beta_{k} \left( q_{ALS}(y_{k-1}) - q_{ALS}(y_{k-2})\right).\label{eq:nesterovALSyonly}
\end{align}
Replacing gradient directions by
update directions provided by ALS is essentially also an approach that has been taken to obtain
nonlinear acceleration of ALS by NGMRES, NCG and LBFGS\cite{sterck2012nonlinear,sterck2015nonlinearly,sterck2018nonlinearly}; in the case
of Nesterov's method the procedure is extremely simple and easy to implement.

In most of what follows, we will focus on the update formula in the $y$ variables, and we will simplify notation for the case of ALS acceleration by defining
\begin{align} 
       \bar{y}_{k} = q_{ALS}(y_{k-1})
\end{align}
and writing the update formula for Nesterov-accelerated ALS as 
\begin{align} 
        y_{k} = \bar{y}_{k} + \beta_{k} \left( \bar{y}_{k} - \bar{y}_{k-1}\right).\label{eq:nesterovbaryonly}
\end{align}
Recently, the application of Nesterov
acceleration to ALS for canonical tensor decomposition was considered by Wang et al.\cite{wang2018accelerating} However, they only tried the vanilla Nesterov technique with a standard Nesterov momentum sequence $\beta_k$ and without restarting or line search
mechanisms %(as for the magenta curve in Fig.\ \ref{fig:intro}), 
and not surprisingly they fail to obtain acceleration of ALS.
The main goal of this paper is to investigate suitable choices for the momentum weights $\beta_{k}$ in Eq.\ (\ref{eq:nesterovbaryonly}), and for restart mechanisms that enable effective convergence patterns when applying the Nesterov-ALS method of Eq.\ (\ref{eq:nesterovbaryonly}).

%-----------------------------------------------------------
\subsection{Canonical Tensor Decomposition.}
%-----------------------------------------------------------
Tensor decomposition has wide applications in machine learning, signal
processing, numerical linear algebra, computer vision, natural language
processing and many other fields\cite{kolda2009tensor}.
This paper focuses on the Canonical Polyadic (CP) decomposition of
tensors\cite{kolda2009tensor}, which is also called the CANDECOMP/PARAFAC decomposition.
CP decomposition approximates a given tensor $T \in \mathbb{R}^{I_{1} \times \ldots \times I_{N}}$
by a low-rank tensor composed of a sum of $R$ rank-one terms,
$\widetilde{T} = \sum_{i=1}^{R} a_{1}^{(i)} \circ \ldots \circ a_{N}^{(i)}$,
where $\circ$ is the vector
outer product. Specifically, defining the factor matrices $A_{i} = (a_{i}^{(1)}, \ldots, a_{i}^{(R)})$,
we minimize the error in the Frobenius norm by
considering the objective function
\begin{align}
	f(A_{1}, \ldots, A_{N})=\frac{1}{2}\norm{T - \sum_{i=1}^{R} a_{1}^{(i)} \circ \ldots \circ a_{N}^{(i)}}_{F}^2.
	\label{eq:cp}
\end{align}
Finding efficient methods for computing tensor decomposition is an active area
of research, but the alternating least squares (ALS) algorithm is still considered one of
the most efficient algorithms for CP decomposition.\cite{acar2011scalable}
\rev{Alternative optimization methods for CP that have been considered in the literature include quasi-Newton methods 
such as NCG and LBFGS\cite{acar2011scalable}, nonlinear least-squares approaches using
Gauss--Newton and Levenberg--Marquardt algorithms\cite{paatero1997weighted,tomasi2006comparison,phan2013low},
and stochastic gradient descent\cite{sidiropoulos2017tensor}.}

ALS finds a CP decomposition in an iterative way.
In each iteration, ALS sequentially updates a block of variables at a time by
minimizing expression (\ref{eq:cp}), while keeping the other blocks fixed:
first $A_{1}$ is updated,
then $A_{2}$,
and so on.
Updating a factor matrix $A_{i}$ is a linear least-squares problem that can be solved in
closed form.
\rev{The ALS update equations for the $A_{i}$ can be derived by considering expressions for the
gradient of $f(A_{1}, \ldots, A_{N})$.
For example, for a tensor of order 3, with factor matrices $A, B, C$ of sizes $I \times R$, $J \times R$, $K \times R$, expressions for the gradient of $f$ are given by Acar et al.\cite{acar2011scalable}
\begin{align}
	\frac{\partial f}{\partial A} = - T_1 (B \odot C) + A \left((B^T B)\ast(C^T C)\right),\label{eq:grad1}\\
	\frac{\partial f}{\partial B} = - T_2 (A \odot C) + B \left((A^T A)\ast(C^T C)\right),\label{eq:grad2}\\
	\frac{\partial f}{\partial C} = - T_3 (A \odot B) + C \left((A^T A)\ast(B^T B)\right), \label{eq:grad3}
\end{align}
where $\odot$ denotes the Khatri-Rao product (the `matching columnwise' Kronecker product\cite{kolda2009tensor,acar2011scalable}), $\ast$ denotes component-wise multiplication, and $T_i$ is the matricized version of tensor $T$ in direction $i$.\cite{kolda2009tensor,acar2011scalable}
For example, in the first equation, $T_1$ is of size $I \times JK$, and $B \odot C$ is of size $JK \times R$.
The $I \times R$ matrix $T_1 (B \odot C)$ is called a `matricized-tensor times Khatri-Rao product' (often abbreviated as
`MTTKRP').}

\rev{The update equations for ALS can be derived from setting each of the gradient components in (\ref{eq:grad1})-(\ref{eq:grad3}) equal to zero\cite{acar2011scalable}:
\begin{align}
	T_1 (B \odot C) = A \left((B^T B)\ast(C^T C)\right),\label{eq:ALS1}\\
	T_2 (A \odot C) = B \left((A^T A)\ast(C^T C)\right),\label{eq:ALS2}\\
	T_3 (A \odot B) = C \left((A^T A)\ast(B^T B)\right).\label{eq:ALS3}
\end{align}
In each ALS iteration, the $A,B,C$ factor matrices are sequentially updated by solving these systems in order for $A, B$ and $C$. It can be shown that these are the normal equations for minimizing $f(A,B,C)$ in $A$, $B$, and $C$, respectively, while keeping the other factor matrices constant\cite{acar2011scalable}. For example, Eq.\ (\ref{eq:ALS1}) are the normal equations for minimizing the least-squares functional $f(A,B,C)=\frac{1}{2}\norm{T - \widetilde{T}}_{F}^2=\frac{1}{2}\norm{T_1-A (B \odot C)^T}_F^2$ with respect to $A$, keeping $B$ and $C$ fixed.}

Collecting the matrix elements of the $A_{i}$'s in a vector $x$,
we use $q_{ALS}(x)$ to denote the updated variables after performing one full ALS 
iteration starting from $x$.

When the CP decomposition problem is ill-conditioned, ALS can be slow to converge
\cite{acar2011scalable}, and
recently a number of methods have been proposed to accelerate ALS.
This includes acceleration by previously mentioned methods such as
NGMRES\cite{sterck2012nonlinear},
NCG\cite{sterck2015nonlinearly}, and LBFGS\cite{sterck2018nonlinearly}.
An approach has also been proposed based on the Aitken-Stefensen acceleration technique.\cite{wang2018accelerating}
These acceleration techniques can substantially improve ALS convergence
speed when problems are ill-conditioned or an accurate solution is required.

%-----------------------------------------------------------
\subsection{Main Approach and Contributions of this Paper.}
%-----------------------------------------------------------

As explained above, our basic approach is to apply Nesterov acceleration to ALS in a manner that is
equivalent to replacing the gradient update in the second step of Nesterov's
method, Eq.\ (\ref{eq:nesterov}), by an ALS step. 
However, applying this procedure directly fails for several reasons.
First, it is not clear to which extent the $\beta_{k}$ momentum weight sequence 
of (\ref{eq:momentum1}), which guarantees optimal convergence
for gradient acceleration in the convex case, applies at all to our case of ALS acceleration
for a non-convex problem. Second, and more generally, it is well-known that optimization
methods for non-convex problems require mechanisms to safeguard against `bad steps',
especially when the solution is not close to a local minimum\cite{nocedal2006numerical}. The main contribution of this
paper is to propose and explore restart-based safeguarding mechanisms for Nesterov acceleration
applied to ALS, along with momentum weight selection.
This leads to a family of acceleration methods for ALS that are competitive with
or outperform previously described highly efficient nonlinear acceleration methods for ALS.
We also compare with choosing the momentum weight using line searches, which is another way
to obtain a safeguarding mechanism and is equivalent to line search methods for ALS that go back to the 1970s \cite{harshman1970foundations,rajih2008enhanced,chen2011new,sorber2016exact}.

\begin{figure}[h!]
	\centering{\includegraphics[width=.7\linewidth]{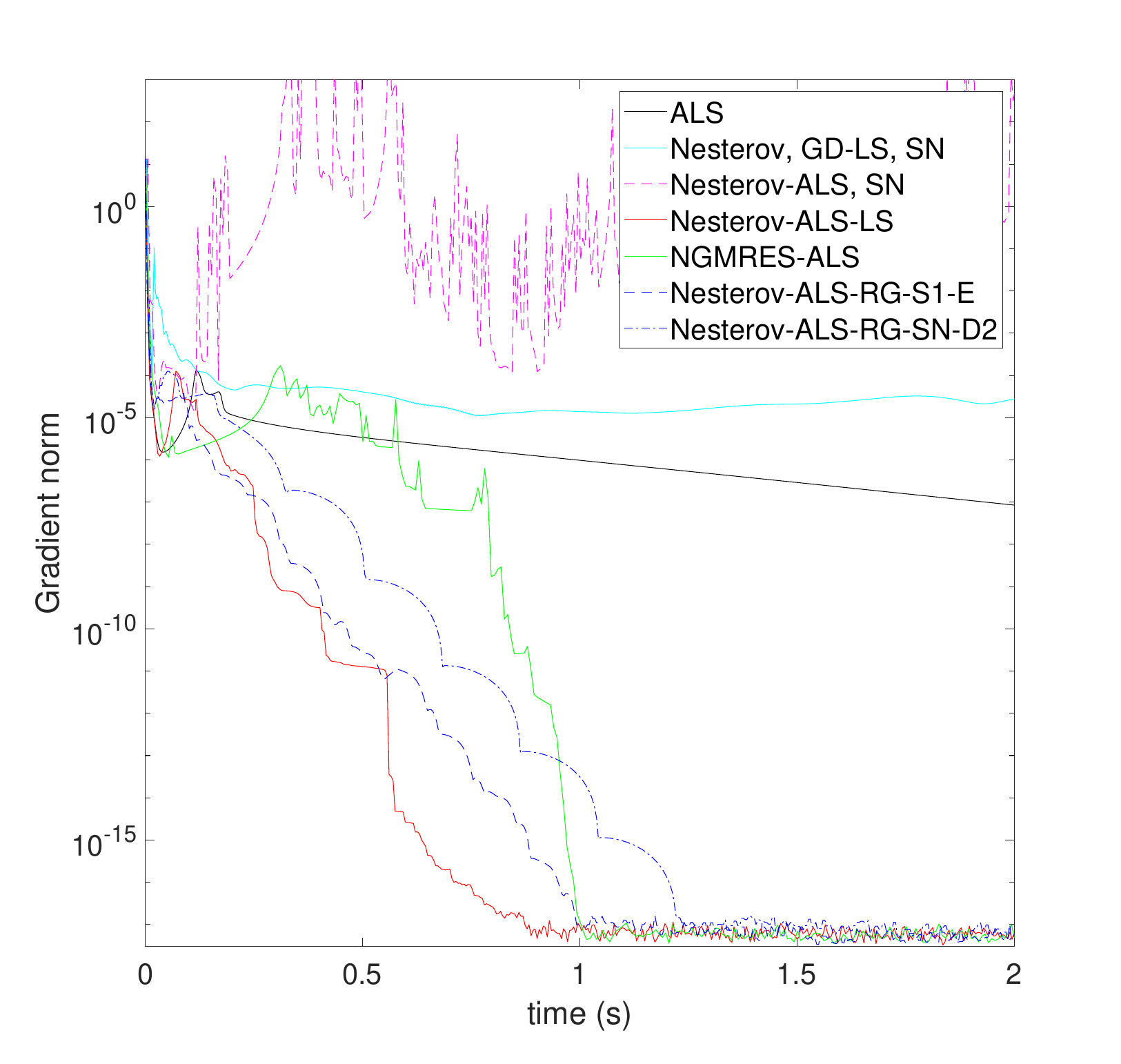}}
	\caption{
Convergence of the gradient norm as a function of runtime for a standard
	ill-conditioned synthetic CP tensor decomposition problem of a tensor of size
	$50 \times 50 \times 50$ (parameters $s$= 50, $c$= 0.9, $R$= 3, $l_1$= 1,
	$l_2$= 1, see \Cref{subsec:synthetic} for the problem description). After an initial transient period, ALS (black) converges at a slow linear rate for this problem, and Nesterov's accelerated gradient method (cyan) converges even more slowly. Direct application of Nesterov acceleration to ALS fails to converge (magenta). To remedy this, we propose a family of methods that relie on restarting mechanisms to stabilize Nesterov-ALS (blue curves), which are competitive with previously proposed nonlinear acceleration methods for ALS, such as NGMRES-ALS (green), but are much easier to implement. We also compare with a Nesterov-ALS method with line search (red). This paper will show numerically that the proposed Nesterov-ALS methods are competitive with or outperform existing nonlinear acceleration methods for ALS, on a set of representative test problems.
	} 
	\label{fig:intro}
\end{figure}

\rev{As further motivation for the problem we address and for our approach,} \Cref{fig:intro} illustrates the convergence difficulties that ALS may experience for ill-conditioned CP tensor decomposition problems, and how nonlinear acceleration may allow to remove these
convergence difficulties. %We address this slow convergence in this paper, and the potential improvements that nonlinear acceleration can provide. 
For the standard ill-conditioned synthetic test problem that is the focus of
Fig.\ \ref{fig:intro} (see \Cref{sec:expt} for the problem description), ALS converges slowly (black curve). 
It is known that standard gradient-based methods such as gradient descent (GD), NCG or LBFGS that do not rely on ALS, perform more poorly than ALS\cite{acar2011scalable}, so it is no surprise that applying Nesterov's accelerated gradient method to the problem (for example, with the gradient descent step length $\alpha_k$ determined by a standard cubic line search\cite{acar2011scalable,sterck2012nonlinear,sterck2015nonlinearly,sterck2018nonlinearly}, cyan curve) leads to worse performance than ALS. 
Nonlinear acceleration of ALS, however, can substantially improve convergence\cite{sterck2012nonlinear,sterck2015nonlinearly,sterck2018nonlinearly}, and we pursue this
using Nesterov acceleration in this paper.
However, as expected, applying Nesterov acceleration (with
the standard Nesterov momentum weight sequence)
directly to ALS for our non-convex problem, by replacing the gradient step
in the Nesterov formula by a step in the ALS direction, does not work and leads to erratic convergence behaviour (magenta curve). 

As the main contribution of this paper, we consider restart mechanisms to stabilize the convergence behaviour of Nesterov-accelerated ALS, and we study how two key parameters, the
momentum step and the restart condition, should be set.
%inspired by previously proposed restarting mechanisms for Nesterov acceleration, extending them to the non-convex and ALS settings. 
The blue curves in Fig.~\ref{fig:intro} show two examples of the acceleration
that can be provided by two variants of the family of restarted Nesterov-ALS
methods we consider. One of these variants (Nesterov-ALS-RG-SN-D2) uses
Nesterov's sequence for the momentum weights, and another successful variant
simply always uses momentum weight one (Nesterov-ALS-RG-S1-E). The naming scheme
for the Nesterov-ALS variants that we consider will be explained in \Cref{sec:expt}. 
Note that, for our non-convex problem, the cascading convergence pattern of our restarted variant with Nesterov's sequence for the momentum weights, Nesterov-ALS-RG-SN-D2, is very similar to the convergence patterns observed in the work of Su et al.\cite{su2016differential}
on restart mechanisms for Nesterov acceleration (using Nesterov's sequence) in the convex setting.
Extensive numerical tests to be provided in \Cref{sec:expt}
show that the best-performing Nesterov-ALS scheme is achieved when using the gradient ratio as momentum weight (as in Nguyen et al.\cite{nguyen2018accelerated}), and restarting when the objective value increases. 
We also compare with determining the Nesterov momentum weight $\beta_k$ in each iteration using a cubic line search (LS) (red curve). The resulting Nesterov-ALS-LS method (which is similar to classical line search methods for ALS \cite{harshman1970foundations,rajih2008enhanced,chen2011new,sorber2016exact}) is competitive with or superior to other recently developed nonlinear acceleration techniques for ALS that use line searches, such as NGMRES-ALS (green curve), with the advantage that Nesterov-ALS-LS is much easier to implement. However, the line searches may require multiple evaluations of $f(x)$ and its gradient and can be expensive. 

The convergence theory of Nesterov's accelerated gradient method for convex problems
does not apply in our case due to the non-convex setting of the CP problem, and because we accelerate ALS steps instead of gradient steps. In fact, in the context of nonlinear convergence acceleration for ALS, few theoretical results on convergence are available \cite{sterck2012nonlinear,sterck2015nonlinearly,sterck2018nonlinearly}.
We will, however, demonstrate numerically, for representative synthetic and real-world
test problems, that our Nesterov-accelerated ALS methods are competitive with or outperform
existing acceleration methods for ALS.
In particular, some of the Nesterov-ALS methods substantially outperform other acceleration methods for ALS when applied to a large real-world ill-conditioned 71$\times$1000$\times$900 tensor.

The remainder of this paper is structured as follows.
%\Cref{sec:related} discusses related work.
\Cref{sec:algo} presents our general Nesterov-ALS scheme and discusses its
instantiations, focusing on the choice of momentum weights and restarting mechanisms.
\Cref{sec:related} discusses and establishes links between the Nesterov-accelerated ALS methods we consider, and existing acceleration methods for ALS.
In \Cref{sec:expt}, we perform an extensive experimental study of our algorithm by
comparing it with a number of acceleration schemes on several benchmark
datasets.
\Cref{sec:conclude} concludes the paper.

%%%%%%%%%%%%%%%%%%%%%%%%%%%%%%%%%%%%%
\section{Nesterov-Accelerated ALS Methods} \label{sec:algo}
%%%%%%%%%%%%%%%%%%%%%%%%%%%%%%%%%%%%%

We consider Nesterov-type acceleration of ALS as in Eq.\ (\ref{eq:nesterovALSyonly}), but a direct
application of a standard Nesterov momentum weight sequence $\beta_k$ for convex problems as in
Eqs.\ (\ref{eq:lambda}--\ref{eq:momentum1}) does not work. 
A typical behavior is illustrated by the magenta curve in \Cref{fig:intro}, which
suggests that the algorithm gets stuck in a highly suboptimal region.
Such erratic behavior arises here, and not in Nesterov's accelerated gradient descent for convex problems, because the ALS update we use for our non-convex problem
is very different from gradient descent.
Below we propose a general restart method to safeguard against bad steps, and investigate suitable choices
for the momentum weights $\beta_k$.

%-----------------------------------------------------------
\subsection{Nesterov-ALS with Restart.}
%-----------------------------------------------------------
Our general Nesterov-ALS scheme with restart is shown in \Cref{alg:nes-als}.
Besides incorporating the momentum term in the update rule (line 12), there are two other
important ingredients in our algorithm:
adaptive restarting (line 5-7), and adaptive momentum weight $\beta_k$ (line 9).
The precise expressions we use for restarting and computing the momentum weight are explained in the following subsections.
In each iteration $k$ of the algorithm we compute a new update according to
the update rule (\ref{eq:nesterovALSxonly}) with momentum term (line 12).
Before computing the update, we check whether a restart is needed (line 5) due
to a bad current iterate. When we restart, we discard the current bad iterate (line 6),
and compute a simple ALS update instead (ALS always reduces $f(x)$ and is thus
well-behaved), by setting $\beta_k$ equal to zero (line 7) such that (line 12) computes
an ALS update. 
Note that, when a bad iterate is discarded, we don't decrease the iteration
index $k$ by one, but instead set the current iterate $x_k$ equal to the previously accepted iterate $x_{k-1}$, which then occurs twice in the sequence of iterates. We wrote the algorithm
down this way because we can then use $k$ to count work (properly keeping track of the cost
to compute the rejected iterate), but the algorithm can of course also be written without duplicating the previous iterate when an iterate is rejected. The index $i$ keeps track of the number of iterates since restarting, which is used for some of our strategies to compute the momentum weight $\beta_k$, see \Cref{sec:momentum}.
The $\beta_{k-1} \neq 0$ condition is required in (line 5), which checks whether a restart is
needed, to make sure that each restart (computing an ALS iteration) is followed by at least
one other iteration before another restart can be triggered (because otherwise the algorithm
could get stuck in the same iterate).

\begin{algorithm}
	\begin{algorithmic}[1]
		\State initialize $x_{1}$, $x_{2} \gets q_{ALS}(x_{1})$
		\State $i \gets 2$ \Comment{$i$ is the number of iterates since \emph{restarting}}
		\State $\beta_{1} =0$ %\Comment{iterate $x_2$ was obtained by ALS, without acceleration}
		\For{$k=2, 3, \ldots \ $} \Comment{$k$ is the number of iterates since the start of the algorithm}

			\If{(restart condition met) \textbf{and} ($\beta_{k-1} \neq 0$)}
				\State $x_{k} = x_{k-1}$ \Comment{discard the current bad step}
				\State $\beta_{k} = 0$, $i \gets 1$ \Comment{force ALS step this iteration}
			\Else
				\State compute $\beta_{k}$ using $i$ and/or previous iterates
			\EndIf
			\State {\textbf{exit} loop if termination criterion is met} 
			\State %$x_{k+1} \gets ALS(x_{k} + \beta_{k} (x_{k} - x_{k-1}))$, $i \gets i + 1$
			$x_{k+1} =q_{ALS}\left(x_{k} + \beta_{k} (x_{k} - x_{k-1})\right)$, $i \gets i + 1$
%			\State {update restart condition parameters if needed} 
% leave this out, because it can be considered part of the "restart condition met" test
		\EndFor
	\end{algorithmic}
	\caption{Nesterov-ALS with restart}
	\label{alg:nes-als}
\end{algorithm}

%We perform a restart whenever the restarting condition is satisfied, and the
%previous momentum weight $\beta_{k-1}$ is not 0.
%The restart sets the momentum weight for the next update to be 0, thus the next
%update is just doing a usual ALS update.
%We require the condition $\beta_{k-1} \neq 0$, because for some problems, the
%restarting condition is always triggered, trapping the algorithm in a loop of doing the same iteration. This occurs for gradient based restarting methods as the ALS algorithm is not guaranteed to reduce the gradient norm, as such we are forced to restart putting us back to the same iterate value as before, we then restart again and adopt the same value with will trigger another restart. As such we become trapped restarting
%If the restarting condition is not satisfied, then we compute an adaptive
%momentum weight.
%This can use the previous iterates, and/or the number of iterates that we have
%seen since last restart.

Various termination criteria may be used. 
In our experiments, we terminate when the gradient 2-norm reaches a set
tolerance:
$$
\norm{\grad f(x_k)} /n_{X} \leq tol, \label{eq:terminate}
$$
Here $n_{X}$ is the number of variables in the low-rank tensor approximation.

The momentum weight $\beta_{k}$ and the restart condition need to be specified
to turn the scheme into concrete algorithms.
We discuss the choices used in \Cref{sec:momentum} and \Cref{sec:restart} below.

%-----------------------------------------------------------
\subsection{Momentum Weight Choices for Nesterov-ALS with Restart.} \label{sec:momentum}
%-----------------------------------------------------------
Naturally, we can ask whether a momentum weight sequence that guarantees optimal convergence for convex problems is applicable in our case.
We consider the momentum weight rule defined in \Cref{eq:momentum1}, but adapted
to take restart into account in \Cref{alg:nes-als}:
\begin{align}
	\beta_{k} &\gets (\lambda_{i-1} - 1)/\lambda_{i}, %& \text{(Fixed weights A)} 
	\label{eq:fixedw}
\end{align}
where $\lambda_{i}$ is defined in \Cref{eq:lambda}.
Restart is taken into account by using $i$ instead of $k$ as the index on the
RHS.

%The second one is defined by
%\begin{align}
%	\beta_{k} &\gets (i-2)/(i+1).  & \text{(Fixed weights B)}
%\end{align}
%\todo{do we used both A and B? if not, comment out the one that is not used. (for Drew)}

%The above momentum weights are in fact not adaptive.
Following Nguyen et al.\cite{nguyen2018accelerated}, we also consider using the 
\emph{gradient ratio} as the momentum weight
\begin{align}
	\beta_{k} &\gets \frac{\norm{\grad f(x_{k})}}{\norm{\grad f(x_{k-1})}}. \label{eq:grw}
\end{align}
%where $f(x)$ is the sum of squared error as defined in \Cref{eq:cp}.
This momentum weight rule can be motivated as follows\cite{nguyen2018accelerated}.
%The momentum weight is large when the algorithm is not making 
When the gradient norm drops significantly, that is, when convergence is fast, the algorithm performs a step closer to the ALS update, because momentum may not really be needed
and may in fact be detrimental, potentially leading to overshoots and oscillations.
When the gradient norm does not change much, that is, when the algorithm is not
making much progress, acceleration may be beneficial and a $\beta_k$ closer to 1 is
obtained by the formula.

Finally, since we observe that Nesterov's sequence \Cref{eq:momentum1}
produces $\beta_k$ values that are always of the order of 1 and approach 1 steadily
as $k$ increases, we can simply consider a choice of $\beta_k=1$ for our non-convex
problems, where we rely on
the restart mechanism to correct any bad iterates that may result, replacing them by
an ALS step. Perhaps surprisingly, the numerical results to be presented below
show that this simplest of choices for $\beta_k$ may work well, if combined with
suitable restart conditions.

%-----------------------------------------------------------
\subsection{Restart Conditions for Nesterov-ALS.} \label{sec:restart}
%-----------------------------------------------------------
One natural restarting strategy is \emph{function restarting} (see, e.g., O'donoghue and Candes\cite{odonoghue2015adaptive,su2016differential} for its use in the convex setting), which restarts
when the algorithm fails to sufficiently decrease the function value.
We consider condition
\begin{align}
	f(x_{k}) &> \eta f(x_{k-d}). \label{eq:fr}
\end{align}
Here, we normally use $d=1$, but $d>1$ can be used to allow for \emph{delay}.
We normally take $\eta=1$, but we have found that it sometimes pays off to allow
for modest increase in $f(x)$ before restarting, and a value of $\eta>1$ facilitates that.
If $d=1$ and $\eta=1$, the condition guarantees that the algorithm will make some
progress in each iteration, because the ALS step that is carried out after a restart
is guaranteed to decrease $f(x)$. However, requiring strict decrease may preclude
accelerated iterates (the first accelerated iterate may always be rejected in favor of an
ALS update), so either $\eta>1$ or $d>1$ allows for a few accelerated iterates to initially
increase $f(x)$, after which they may decrease $f(x)$ in further iterations in a much
faster way than ALS, potentially resulting in substantial acceleration of ALS.
%\todo{Nan, can you add references for the statement in the previous sentence?}
While function restarting (with $d=1$ and $\eta=1$) has been observed to significantly improve convergence for convex problems, no theoretical convergence rate has been obtained\cite{odonoghue2015adaptive,su2016differential}.

Following Su et al.\cite{su2016differential}, we also consider the \emph{speed
restarting} strategy which restarts when 
\begin{align}
	\norm{x_{k} - x_{k-1}} &<  \norm{x_{k-1} - x_{k-2}}. \label{eq:sr}
\end{align}
Intuitively, this condition means that the speed along the convergence trajectory,
as measured by the change in $x$, drops.
Su et al.\cite{su2016differential} showed that speed restarting leads to guaranteed
improvement in convergence rate for convex problems.

Another natural strategy is to restart when the gradient norm satisfies
\begin{align}
	\norm{\grad f(x_{k})} &> \eta \norm{\grad f(x_{k-d})}, \label{eq:gr}
\end{align}
where, as above, $\eta$ can be chosen to be equal to or greater than one.
This \emph{gradient restarting} strategy (with $\eta=1$) has been used in conjunction with
gradient ratio momentum weight by Nguyen et al.\cite{nguyen2018accelerated}, and
a similar condition on the residual has been used for Nesterov acceleration of ADMM
for convex problems by Goldstein et al.\cite{goldstein2014fast}.

When we use a value of $\eta>1$ in the above restart conditions,
we have found in our experiments that it pays off
to allow for a larger $\eta$ immediately after the restart, and then decrease $\eta$ in subsequent steps. 
In particular, in our numerical tests below, we set $\eta$ = 1.25, and decrease $\eta$ in
every subsequent step by 0.02, until $\eta$ reaches 1.15.

%-----------------------------------------------------------
\subsection{Nesterov-ALS with Line Search.}
%-----------------------------------------------------------
To compare our numerical results with existing acceleration methods for ALS that use a line search in the direction
of the ALS update\cite{harshman1970foundations,rajih2008enhanced,chen2011new,sterck2012nonlinear,sorber2016exact}, we also consider an approach where the momentum weight $\beta_k$ in the Nesterov extrapolation formula (\ref{eq:nesterovALSxonly}) is determined by a line search.
In Section \ref{subsec:compare-line-search} we explain the equivalence of this approach with existing line search methods to
accelerate ALS\cite{harshman1970foundations,rajih2008enhanced,chen2011new,sterck2012nonlinear,sorber2016exact}.
The line search to determine the momentum weight $\beta_k$ safeguards against
bad steps introduced by the $\beta_{k} (x_{k} - x_{k-1})$ term, so an additional restart mechanism is not needed.
(Note that ALS itself always reduces $f(x)$ and is not prone to introducing bad steps; if the line search does not find a suitable extrapolation point, it simply returns $\beta_{k}=0$ and the ALS step is accepted, which can be considered as a restart mechanism
built-in into the line search.)
For the line search approach, we determine $\beta_{k}$ in each iteration as an approximate solution of
\begin{align}
	\beta_{k} & \approx \arg\min_{\beta \ge 0} f(x_{k} + \beta (x_{k} - x_{k-1})). %\\
	%&\qquad\qquad \text{(Line search)} \nonumber
	\label{eq:lsw}
\end{align}
%This provides an optimal greedy method to choose the momentum weight.
%Indeed, our experiments show that it can significantly reduce the number of
%iterations required for convergence.
We use the standard Mor\'{e}-Thuente cubic line search that has been used extensively 
for tensor decomposition methods\cite{acar2011scalable,sterck2012nonlinear,sterck2015nonlinearly,sterck2018nonlinearly}.
This inexact line search finds a value of $\beta_k$ that satisfies the Wolfe conditions, which impose a sufficient descent condition and a curvature condition. Each iteration of this iterative
line search requires the computation of the function value, $f(x)$, and its gradient. As such, the line search can be quite expensive. In our numerical tests, we use the following line search parameters: $10^{-4}$ for the descent condition, $10^{-2}$ for the curvature condition, a starting search step length of 1, and a maximum of 20 line search iterations.

%In particular, in our numerical tests below, we set $\eta=1.15$, and
%decrease $\eta$ in every subsequent step by 0.02, until $\eta$ reaches 1.
%\todo{Drew, is the previous sentence correct?}

%Another gradient-based restarting strategy is to restart when
%\begin{align}
	%\grad f(x_{k-1} + \beta_{k-1} (x_{k-1} - x_{k-2}))^{\top} (x_{k} - x_{k-1}) < 0.
%\end{align}
%This is adapted from the restart criterion 
%$\grad f(y_{k-1})^{\top} (x_{k}) - x_{k-1}) < 0$ used for the accelerated gradient method
%as in \cite{odonoghue2015adaptive}.

%\comment{\medskip\noindent{\bf Continuous time dynamics.} 
%We derive a continuous time dynamics of Nesterov accelerated ALS by taking
%small steps towards the next iterate obtained using the ALS operator, that is,
%we consider the following Nesterov accelerated ALS
%\begin{align}
%	x_{k+1} = ALS_{s}(x_{k} + \beta_{k} (x_{k} - x_{k-1})).
%\end{align}
%where $ALS_{s} = (1 - s) \I + s ALS$ with $\I$ being the identity operator.
%The update $x' = ALS_{s}(x)$ can be written as $x' = x - s(x - ALS(x))$, that
%is, it moves along the direction $x - ALS(x)$ for a step size of $s$.
%
%Let $X(t) = X(k\sqrt{s}) = x_{k}$, and $\beta_{k} = \frac{k-2}{k+1}$, then
%following a similar derivation in \cite{su2016differential}, as $s \to 0$, 
%we have
%\begin{align}
%	\ddot{X} + \frac{3}{t} \dot{X} + Y - ALS(X) = 0.
%\end{align}
%
%\todo{Plot function value against $k\sqrt{s}$ for this variant of accelerated ALS using small
%$s$.
%Check whether this exhibits oscillatory behavior.
%If yes, then generate the plot when restart is used.
%}}

%%%%%%%%%%%%%%%%%%%%%%%%%%%%%%%%%%%%%%%%%%%%%%
\section{Relation with Existing Acceleration Methods} \label{sec:related}
%%%%%%%%%%%%%%%%%%%%%%%%%%%%%%%%%%%%%%%%%%%%%%

In this section we elucidate the relation of the Nesterov acceleration methods for ALS discussed in this paper with other existing convergence acceleration methods for ALS.

%-----------------------------------------------------------
\subsection{NGMRES and Anderson Acceleration}
%-----------------------------------------------------------
The Nonlinear Generalized Minimal Residual Method\cite{sterck2012nonlinear,sterck2013steepest} (NGMRES) for minimizing $f(y)$
accelerates convergence of a nonlinear iteration $y_k=q(y_{k-1})$ by considering the $w$-step extrapolation formula
\begin{align} 
        \widehat{y}_{k} = \bar{y}_{k} + \sum_{j=k-w}^{k-1} \beta_j (\bar{y}_k-y_j) .\label{eq:NGMRES}
\end{align}
where $\bar{y}_k=q(y_{k-1})$,
in combination with a line search to stabilize the iteration when the iterate is far from a stationary point,
\begin{align} 
        y_k = \bar{y}_k + \alpha_k (\widehat{y}_k-\bar{y}_k) ,\label{eq:NGMRESLS}
\end{align}
where $\alpha_k$ is determined by a line search.
The expansion coefficients $\beta_j$ are computed in each step by solving a small $w \times w$ linear least-squares problem
that minimizes a linearization of $\|\nabla f(\widehat{y}_k)\|_2$.
NGMRES was originally proposed as a convergence accelerator for solving a system of nonlinear algebraic equations $g(y)=0$, and it was shown to be essentially equivalent to the GMRES method in the case of a linear system $Ay=b$\cite{washio1997krylov,oosterlee2000krylov}.
\rev{Note that the nonlinear iteration $y_k=q(y_{k-1})$ can also be seen as an inner-iteration nonlinear preconditioner for NGMRES, see also Section  \ref{subsec:nonl_prec} below.}

NGMRES is closely related to the classical Anderson acceleration method for solving $g(y)=0$,\cite{brune2015composing} which uses the expansion formula
\begin{align} 
        \widehat{y}_{k} = \bar{y}_{k} + \sum_{j=k-w}^{k-1} \beta_j (\bar{y}_k-\bar{y}_j) \label{eq:Anderson},
\end{align}
where the $\beta_j$ are also determined to approximately minimize $\|g(\widehat{y}_k)\|_2$. Anderson 
acceleration can be viewed as a multi-secant method\cite{fang2009two} and also reduces to
GMRES for the case of linear systems\cite{walker2011anderson}.

It can be seen immediately that the from of the Nesterov extrapolation formula (\ref{eq:nesterovbaryonly}) is equivalent to the Anderson extrapolation formula (\ref{eq:Anderson}) with one step ($w=1$), and is also closely related to NGMRES with one step (Eq.\ (\ref{eq:NGMRES})). The
difference, of course, is that the extrapolation weight $\beta_k$ in Nesterov's formula (\ref{eq:nesterovbaryonly}) is usually determined differently than for Anderson acceleration and NGMRES.

%-----------------------------------------------------------
\subsection{Line Search Methods for ALS}\label{subsec:compare-line-search}
%-----------------------------------------------------------
As early as 1970, Harshman\cite{harshman1970foundations} described enhancing ALS convergence for CP tensor decomposition by an over-relaxation in the direction of the ALS update:
\begin{align} 
y_k=\bar{y}_k+ \beta (\bar{y}_k-y_{k-1}) \label{eq:harshman-overrelaxation},
\end{align}
where $\bar{y}_k=q_{ALS}(y_{k-1})$ and a value between 1.2 and 1.3 was recommended for the over-relaxation parameter $\beta$.
This extrapolation step is of the same form as NGMRES with window size one, see Eq.\ (\ref{eq:NGMRES}), and is similar to (but not the same as) Nesterov extrapolation formula (\ref{eq:nesterovbaryonly}) and Anderson acceleration with $w=1$.
In later work, $\beta$ in Eq.\ (\ref{eq:harshman-overrelaxation}) was determined in each step using a line search, obtaining a value of $\beta$ that approximately minimizes $f(y_k)$. 
In particular, the so-called \emph{enhanced line search} approach\cite{rajih2008enhanced} computes the optimal $\beta$ in Eq.\ (\ref{eq:harshman-overrelaxation}) in an accurate manner, which is possible because, for finding the optimal CP decomposition in the 2-norm, the objective is a polynomial function of $\beta$.
Enhanced line search was later extended to broader tensor optimization problems and to exact plane search\cite{sorber2016exact}.
Chen et al.\cite{chen2011new} consider one-step extrapolation according to the update formula
\begin{align} 
y_k=\bar{y}_k+ \beta_k (\bar{y}_k-\bar{y}_{k-1}), \label{eq:chen-1step}
\end{align}
where again $\bar{y}_k=q_{ALS}(y_{k-1})$, which is of the same form as the Nesterov-type update we consider in this paper, Eq.\ (\ref{eq:nesterovbaryonly}) (and, thus, the same as Anderson acceleration with one step), but Chen et al.\cite{chen2011new} determine $\beta_k$ using a line search. Chen et al.\cite{chen2011new} also considered two variants of two-step extrapolation formula
\begin{align} 
y_k=\bar{y}_k+ \beta_k (\gamma_0 \bar{y}_k+\gamma_1 \bar{y}_{k-1}+\gamma_2 \bar{y}_{k-2}), \label{eq:chen-2step},
\end{align}
where the coefficients $\gamma_j$ in the extrapolation are pre-determined based on a so-called geometric or algebraic ansatz. This approach is similar to NGMRES or Anderson acceleration with window size two, except that the latter methods aim to determine optimal extrapolation coefficients by solving, in each iteration, a $w \times w$ linear least-squares problem for a window size of $w$.
NGMRES and Anderson, thus, adapt the extrapolation coefficients each iteration to minimize the residual of the gradient, for general window size, whereas the methods of Chen et al.\cite{chen2011new} are limited to two-step extrapolation and use extrapolation coefficients that are fixed over the iterations and are determined in an ad-hoc way.
Numerically it has been shown\cite{sterck2012nonlinear} that NGMRES with window size greater than one (i.e., multistep extrapolation) can be much faster than (enhanced) line search of the form (\ref{eq:harshman-overrelaxation}) (i.e., one-step extrapolation); it was found that NGMRES window size $w=20$ is a good choice for efficient acceleration of ALS for CP tensor decomposition\cite{sterck2012nonlinear}.
Finally, we note that acceleration of ALS by NGMRES, NCG and LBFGS also employs line searches as a globalization mechanism to guard against erratic iterates\cite{sterck2012nonlinear,sterck2015nonlinearly,sterck2018nonlinearly}.

While the main focus of this paper is on investigating momentum weights $\beta_k$ and restart mechanisms for Nesterov update formula (\ref{eq:nesterovbaryonly}) for ALS, we also compare with a version of update formula (\ref{eq:nesterovbaryonly}) where $\beta_k$ is determined by a line search, similar to existing methods for accelerating ALS convergence that use line searches
\cite{harshman1970foundations,rajih2008enhanced,chen2011new,sterck2012nonlinear,sorber2016exact}.

%-----------------------------------------------------------
\subsection{Nonlinear Preconditioning: Left Preconditioning and Right Preconditioning} \label{subsec:nonl_prec}
%-----------------------------------------------------------
Applying Nesterov's acceleration approach to ALS as explained in Section \ref{subsec:Nesterov-accel} can also be interpreted in the context of nonlinear preconditioning\cite{brune2015composing} for optimization methods. In recent work\cite{sterck2018nonlinearly} a general framework for nonlinear preconditioning of optimization methods was formulated and applied to nonlinear preconditioning of NCG and LBFGS. 
This framework extends the concepts of linear preconditioning for linear systems
$$
Ay=b
$$
to genuinely nonlinear preconditioners for nonlinear optimization.
In the context of linear systems, the concept of \emph{left preconditioning} multiplies the linear system with a nonsingular preconditioning matrix $P$, aiming to improve the convergence of iterative methods (such as GMRES) applied to the preconditioned system 
$$
PAy=Pb.
$$
Left preconditioning can be seen as taking linear combinations of the equations to be solved.
In the \emph{right preconditioning} approach, on the other hand, a linear change of the variables $y=Pz$ is considered by means of a nonsingular preconditioning matrix $P$, and the iterative method is then applied to the preconditioned system
$$
APz=b.
$$
We now explain how these ideas can be extended to minimizing nonlinear functions $f(y)$ using genuinely nonlinear preconditioners, where the preconditioner is not given by a linear coordinate transformation that is encoded by a matrix multiplication, but by a genuinely nonlinear transformation.

\subsubsection{Nonlinear Left Preconditioning}
\rev{In the linear case, left preconditioning of, for example, the GMRES method, can be understood as a combination of two methods: an inner preconditioning iteration
%..................................................
\begin{align}
y_{k+1}=y_k - P (A y_k-b),
\label{eq:prec_lin}
\end{align}
%..................................................
where $P$ is the left preconditioning matrix,
is combined with the outer GMRES iteration (where each iteration of the combined method typically uses one inner update and one outer update)\cite{washio1997krylov,oosterlee2000krylov,sterck2012nonlinear}.
One can say that inner iteration (\ref{eq:prec_lin}) preconditions the GMRES outer iteration, speeding up the convergence of GMRES, or, in an alternative view, one can say that GMRES is an outer iteration that accelerates the convergence of inner iteration (\ref{eq:prec_lin}). Possible choices for preconditioning matrix $P$ include, for example, the lower triangular part of $A$ (Gauss-Seidel), or the constant diagonal preconditioner $\alpha I$ (Richardson iteration), with a suitably chosen $\alpha$.}

\rev{It is instructive to consider the case of solving a linear system $Ay=b$ where $A$ is Symmetric Positive Definite (SPD). In this case, solving $Ay=b$ is equivalent to minimizing $f(y)=\frac{1}{2}y^TAy-b^Ty$. The gradient $g(y)$ of $f(y)$ is given by
%..................................................
\begin{align}
g(y)= \nabla f(y)=Ay-b.
\label{eq:grad_lin}
\end{align}
%..................................................
In this case, it can be seen that inner iteration preconditioner (\ref{eq:prec_lin}) using Richardson iteration is equivalent to steepest descent:
%..................................................
\begin{align}
y_{k+1}=y_k - \alpha \nabla f(y_k).
\label{eq:steep1}
\end{align}
%..................................................
More generally, when using a preconditioning matrix $P$ in the case of an SPD system, the expression $P (A y_k-b)= P \nabla f(y_k) =P g(y_k)$ in (\ref{eq:prec_lin}) can be called the preconditioned gradient. Note that the preconditioned gradient can be obtained from the preconditioning iteration (\ref{eq:prec_lin}) by
%..................................................
\begin{align}
P g(y_k)=-(y_{k+1}-y_k).
\label{eq:precond_grad_lin}
\end{align}
%..................................................
In words, the preconditioned gradient vector is given by the update vector of the linear preconditioning iteration.}

\rev{Conceptually, this can be generalized to nonlinear preconditioning as follows.
Standard nonlinear optimization methods such as NCG or LBFGS use gradient directions to minimize a nonlinear function $f(y)$.
However, if a simple nonlinear iterative method for minimizing $f(y)$ is available that converges faster than gradient descent, and, in a sense, provides better search directions than the gradient, we can use this iterative method as an inner preconditioning iteration for NCG or LBFGS, and make use of the more suitable search directions provided by this method.
For example, in the case of ALS for canonical tensor decomposition, we can use the nonlinear iteration
%..................................................
\begin{align}
y_{k+1}=q_{ALS}(y_k),
\label{eq:prec_nonlin}
\end{align}
%..................................................
as the inner iteration preconditioner of NCG or LBFGS. In analogy with the linear SPD case, we can interpret the direction provided by ALS as a nonlinearly preconditioned gradient direction, $\mathcal{P}(g;y)$, given by
%..................................................
\begin{align}
\mathcal{P}(g;y_k)= - (y_{k+1}-y_k) = -(q_{ALS}(y_{k})-y_k),
\label{eq:precond_grad_nonlin}
\end{align}
%..................................................
in analogy with (\ref{eq:precond_grad_lin}) and (\ref{eq:prec_lin}). To bring this idea into practice, we simply replace the gradient direction in the update formulas for NCG or LBFGS by the nolinearly preconditioned gradient $\mathcal{P}(g;y)$, i.e., by the update vector that the nonlinear preconditioner (used as an inner iteration) provides. In a sense, standard NCG and LBFGS can be understood as using steepest descent (gradient directions) as the inner iteration. In left-nonlinearly preconditioned NCG and LBFGS, we replace the steepest-descent inner iteration by the nonlinear preconditioner iteration (e.g., ALS). Another way to understand this is that standard NCG and LBFGS are iterative methods that work on solving $g(y)=0$, i.e., their update formulas drive the gradient to zero (equivalent to solving $Ay-b=0$ in the linear case). Nonlinearly preconditioned NCG or LBFGS work on solving the nonlinearly preconditioned equation $\mathcal{P}(g;y)=0$ (which is really the fixed-point equation $y-q(y)=0$), i.e., their update formulas drive the nonlinearly preconditioned gradient to zero (equivalent to solving $PAy-Pb=0$ in the linear case). Further details and numerical results can be found in a recent paper establishing this general formalism of nonlinear left preconditioning for optimization\cite{sterck2018nonlinearly}. In this formalism, we can say that the nonlinear preconditioner is used as an inner iteration for the outer-iteration NCG and LBFGS methods, or, alternatively, we can say that NCG and LBFGS are used as nonlinear convergence accelerators for the inner iteration (e.g., ALS).}

\rev{This nonlinear left preconditioning can in principle be applied to any nonlinear optimization method (not just NCG or LBFGS), and it can clearly also be applied to Nesterov's accelerated gradient descent method as in Eq.\ (\ref{eq:nesterovyonly}), when, for example, using ALS as the nonlinear preconditioner for CP tensor decomposition. In fact, in the case of Nesterov's method, the nonlinear left preconditioning procedure is very simple:
it basically amounts to replacing the gradient update $-\alpha_{k} \grad f(y_{k})$ in \cref{eq:nesterov} by $q_{ALS}(y_k)-y_k$, the step provided by ALS, and directly results in the nonlinearly preconditioned formulas of Eq.\ (\ref{eq:nesterovqxonly}) or Eq.\ (\ref{eq:nesterovqyonly}), so the Nesterov acceleration formula (\ref{eq:nesterovALSyonly}) for ALS can also be interpreted as using ALS as a nonlinear preconditioner for Nesterov's accelerated gradient formula (\ref{eq:nesterovyonly}).}

\subsubsection{Nonlinear Right Preconditioning -- Transformation Preconditioning for Optimization}
\rev{Similarly, nonlinear preconditioning techniques can be derived for optimization that rely on a nonlinear change of variables\cite{sterck2018nonlinearly}, inspired by right preconditioning in the linear case. In the so-called transformation preconditioning approach\cite{sterck2018nonlinearly} (which is a form of right preconditioning),
nonlinearly preconditioned versions of NCG and LBFGS are derived as follows. One first considers minimization of $f(y)$ using a linear change of variables
$$y=Cz,$$
and then applies standard NCG and LBFGS (in the $z$ variable) to minimize
$$\widehat{f}(z)=f(Cz).$$
Transforming the resulting update formulas back to the original $y$ variables using
%..................................................
\begin{align}
\nabla_z \widehat{f}(z) =\nabla_z f(Cz)=C^T \nabla_y f(y)
\label{eq:ztransformgrad}
\end{align}
%..................................................
introduces the linear preconditioning matrix
%..................................................
\begin{align}
P=CC^T
\label{eq:Pcct}
\end{align}
%..................................................
in expressions like the scalar product of gradients:
%..................................................
\begin{align}
\nabla_z \widehat{f}(z)^T \nabla_z \widehat{f}(z)= \nabla_y f(y)^T CC^T \nabla_y f(y)=\nabla_y f(y)^T P \nabla_y f(y).
\label{eq:scalargrad}
\end{align}
%..................................................
The resulting linearly preconditioned NCG and LBFGS methods are well-known\cite{hager2006survey,luenberger1984linear}.
However, in recent work this approach was extended to nonlinear preconditioning by replacing the linearly preconditioned gradient $P\nabla_y f(y_k)$
in expressions like (\ref{eq:scalargrad}) by the nonlinearly preconditioned gradient direction $\mathcal{P}(g;y_k)$ given by the update vector that is provided by the nonlinear preconditioning iteration (\ref{eq:prec_nonlin}): $\mathcal{P}(g;y_k)= - (y_{k+1}-y_k) = -(q_{ALS}(y_{k})-y_k),$ as in (\ref{eq:precond_grad_lin}).
This gives different nonlinearly preconditioned iterations for NCG and LBFGS than the nonlinear left preconditioning approach discussed above.
An attractive feature of this nonlinear transformation preconditioning is that it reduces to well-known linear preconditioning techniques for NCG and LBFGS in the case of the linear change of variables $y=Cz$\cite{hager2006survey,luenberger1984linear}. In practice, both left and transformation nonlinear preconditioning may lead to dramatic improvements in convergence for NCG and LBFGS\cite{sterck2018nonlinearly}. In this paper, we use the transformation preconditioning versions of NCG-ALS and LBFGS-ALS in the numerical results.
In the case of Nesterov's method, which is a simple iteration that does not involve scalar products of gradients, the two procedures of nonlinear left preconditioning and transformation preconditioning give the same result, i.e., both approaches lead to the nonlinearly preconditioned formulas of Eq.\ (\ref{eq:nesterovqxonly}) or Eq.\ (\ref{eq:nesterovqyonly}).}

More broadly, nonlinear preconditioning has a long but not widely known history in computational science, and is currently an active area of ongoing research\cite{brune2015composing}, including in the optimization context\cite{sterck2018nonlinearly}.
In our numerical results, we compare Nesterov acceleration of ALS --- i.e., Nesterov's method with ALS as nonlinear preconditioner ---, with NCG and LBFGS acceleration of ALS --- i.e., NCG and LBFGS with ALS as nonlinear preconditioner.
Numerical results will show that Nesterov-ALS is often competitive with the more sophisticated LBFGS-ALS. In addition, Nesterov acceleration is attractive because it is much easier to implement than LBFGS acceleration.

%%%%%%%%%%%%%%%%%%%%%%%%%%%%%%%%%%%%%%%%%%%%%%
%%%%%%%%%%%%%%%%%%%%%%%%%%%%%%%%%%%%%%%%%%%%%%
\section{Numerical Tests} \label{sec:expt}
%%%%%%%%%%%%%%%%%%%%%%%%%%%%%%%%%%%%%

We evaluated our Nesterov-ALS algorithm with various choices for momentum weights and restart mechanisms on a set of synthetic CP test problems that have been carefully designed and used in many papers, and three real-world datasets of different
sizes and originating from different applications. All numerical tests were performed in
Matlab, using the Tensor Toolbox\cite{TTB_Software} and the Poblano Toolbox for optimization\cite{SAND2010-1422}. 
%Matlab code for our methods and tests will be made available on the authors' webpages, and as an extension to the Poblano Toolbox.

\subsection{Naming Convention for Nesterov-ALS Schemes.}
We use the following naming conventions for
the restarting strategies and momentum weight
strategies defined in \Cref{sec:algo}.
For the restarted Nesterov-ALS schemes, we append Nesterov-ALS with
the abbreviations in \Cref{tbl:abbr} to denote the restarting scheme used, and the
choice for the momentum weight $\beta_k$. 

For example, Nesterov-ALS-RF-SG means using restarting based on function value (RF)
and momentum step based on gradient ratio (SG).
For most tests we don't use delay (i.e., $d=1$ in \Cref{eq:fr} or \Cref{eq:gr}),
and $\eta$ is usually set to 1 in \Cref{eq:fr} or \Cref{eq:gr}.
Appending D$n$ or E to the name indicates that a delay $d=n>1$ is used, and that $\eta\ne1$ is
used, respectively.
The line search Nesterov-ALS scheme we compare with is denoted Nesterov-ALS-LS.

%To indicate the restarting strategy, we use RF, RX, RG to denote the function
%restarting strategy of \Cref{eq:fr}, the speed restarting strategy of
%\Cref{eq:sr} that is based on $x$-updates, or the gradient restarting strategy of \Cref{eq:gr},
%respectively.
%For the momentum step choice $\beta_k$, we use SN, SG, and S1 to denote the use of 
%the fixed Nesterov sequence in \Cref{eq:fixedw}, the gradient ratio in \Cref{eq:gr}, or a
%step of constant size 1, respectively.
%Appending D$n$ or E to the name indicates that the delay $d=n$, and that $\eta\ne1$ is
%used, respectively.

\begin{table}
	\centering
	\begin{tabular}{ll}
		\toprule
		Abbreviation & Explanation \\
		\midrule
		RF & function restarting as in \Cref{eq:fr}\\
		RG & gradient restarting \rev{as in \Cref{eq:gr}} \\
		RX & speed restarting as in \Cref{eq:sr}\\
		\midrule
		SN & Nesterov step as in \Cref{eq:momentum1} \\
		SG & gradient ratio step \rev{as in \Cref{eq:grw}}\\
		S1 & constant step 1\\
		\midrule
		\rev{D$d$} & \rev{delay $d>1$ in restart condition} \\
		E & $\eta\ne1$ in restart condition \\
		\bottomrule
	\end{tabular}
	\caption{Abbreviations used in naming convention for restarted Nesterov-ALS variants.}
	\label{tbl:abbr}
\end{table}

\subsection{Baseline Algorithms.}
We compare our proposed Nesterov-ALS schemes with the recently proposed nonlinear
acceleration methods for ALS using GMRES\cite{sterck2012nonlinear},
NCG\cite{sterck2015nonlinearly}, and LBFGS\cite{sterck2018nonlinearly}.
These methods will be denoted in the result figures as GMRES-ALS, NCG-ALS,
and LBFGS-ALS, respectively.

%\subsection{Implementation details.}
%For each test we randomly generate a seeded initial guess. For the synthetic dataset each test problem has ten different seeded initial guesses. For each test problem we generate a tensor, which are different based on the seed, with the same parameters. For the real world problems we have a set tensor and simply change the initial guess. In these tests the number of different seeds changes. 

\subsection{Synthetic Test Problems and Results.}
\label{subsec:synthetic}

We use the synthetic tensor test problems considered by\cite{acar2011scalable} and
used in many papers as a standard benchmark test problem for CP decomposition\cite{sterck2012nonlinear,sterck2015nonlinearly,sterck2018nonlinearly}.
As described in more detail elsewhere\cite{sterck2012nonlinear}, we generate six classes
of random three-way tensors with highly collinear columns in the factor matrices.
We add two types of random noise to the tensors generated from the factor matrices
(homoscedastic and heteroscedastic noise\cite{acar2011scalable,sterck2012nonlinear}),
and then compute low-rank CP decompositions of the resulting tensors.

Due to the high collinearity, the problems are ill-conditioned and ALS
is slow to converge\cite{acar2011scalable}. All tensors have equal size
$s=I_1=I_2=I_3$ in the three tensor dimensions. The six classes differ
in their choice of tensor sizes ($s=20, \ 50, \ 100$), decomposition rank
($R=3, \ 5$), and noise parameters $l_1$ and $l_2$ ($l_1=0,\ 1$ and
$l_2=0,\ 1$), in combinations that are specified in Table \ref{table:test_problems}
of Appendix \ref{ap:params}.

To compare how various methods perform on these synthetic problems, we generate
10 random tensor instances with an associated random initial guess for each of the
six problem classes, and run each method on each of the 60 test problems,
with a convergence tolerance $tol=10^{-9}$. 
We then present so-called $\tau$-plot performance profiles\cite{sterck2012nonlinear}, as explained below, to compare the relative performance of the methods over the test problem set.

%In our experiment, we generated 12 tensors using different parameters (see the
%appendix for the list of parameters used).
%Six of them have a low collinearity parameter of $c=0.5$ and six of them have a
%high collinearity parameter of $c=0.9$.
%Problems with high collinearity require a much larger number of ALS iterions to
%converge and are considered as ill-conditioned, problems with low collinearity
%are considered well conditioned \cite{acar2011scalable}. 
%For this reason we refer to problems with low collinearity as easy and problems
%with high collinearity as hard. The rank of each decomposition is either 3 or 5. For each size $s$ and collinearity $c$ we attempt both a rank-3 and a rank-5  CP decomposition. See \ref{table: test_problems} for the full list of test problems and corresponding parameter values.
%\todo{explanation on collinearity
%and noise parameters, in the appendix (for Drew)... I'm slightly confused about this, did we want more than the explanation on collinearity as above? Or the same thing in the appendix. Also Hans would you be able to cover the noise? I'm slightly uncertain about how this actually works.}

\medskip\noindent\emph{Optimal restarted Nesterov-ALS method}.
Our extensive experiments on both the synthetic and real-world datasets (as
indicated in tests below and in Appendix \ref{ap:detailed-results}) suggest that the optimal
restarted Nesterov-ALS method is the one using function restarting (RF) and
gradient ratio momentum steps (SG), i.e., Nesterov-ALS-RF-SG.
As a comparison, the study of Nguyen et al.\cite{nguyen2018accelerated}
suggests that gradient restarting and gradient
ratio momentum weights work well for accelerating gradient descent
in the context of nonlinear system solving.

\begin{figure}
%	\begin{subfigure}{\linewidth}
%		\caption{time}
%	\end{subfigure}
%	\begin{subfigure}{\linewidth}
%		\caption{iteration}
%	\end{subfigure}
%	\begin{subfigure}{\linewidth}
		\centering{\includegraphics[width=.7\linewidth]{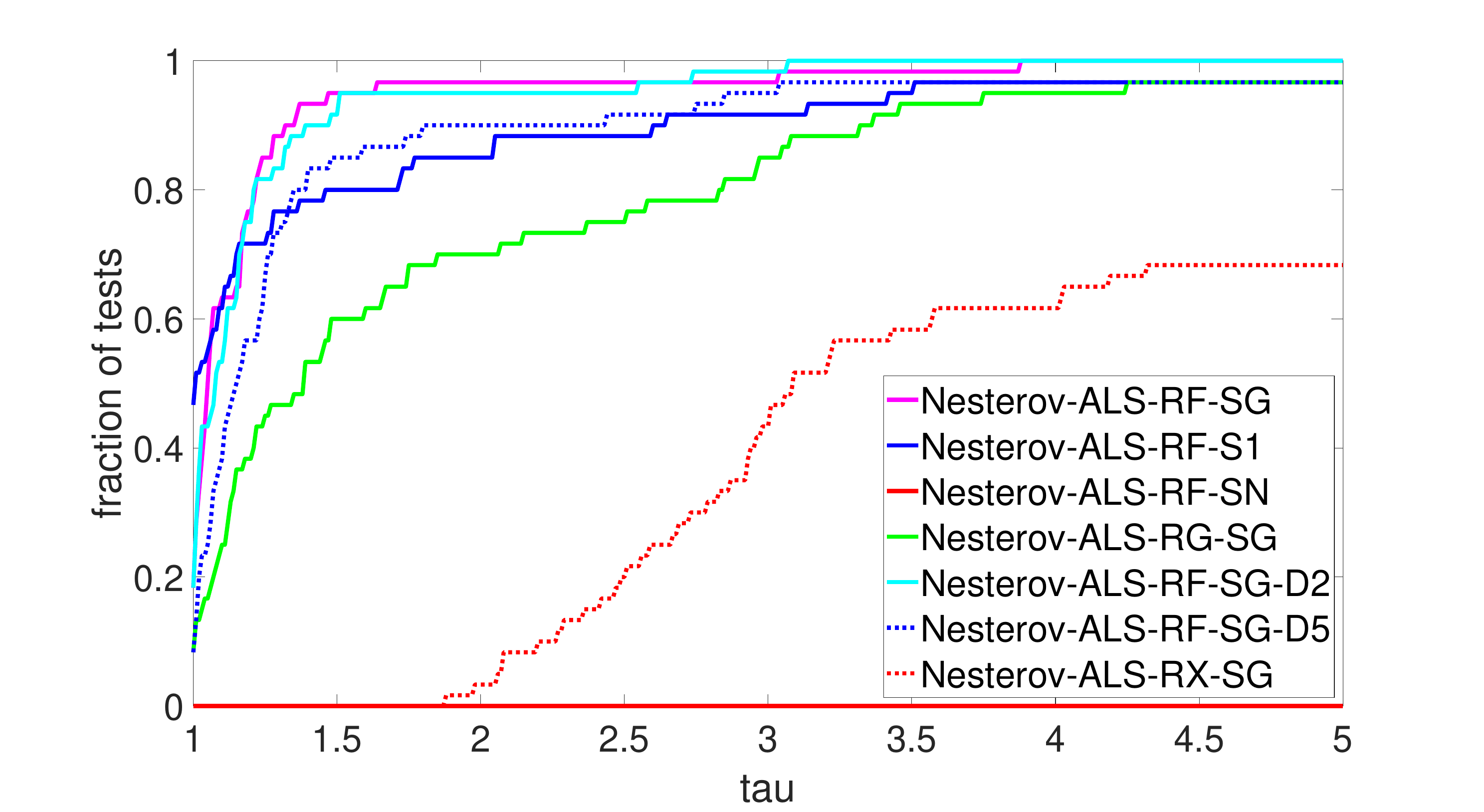}}
%		\caption{tau plot}
%	\end{subfigure}
	\caption{Synthetic test problems. $\tau$-plot comparing the optimal restarted algorithm, Nesterov-ALS-RF-SG, with other variants of restarted Nesterov-ALS.}
	\label{fig:compare-nes-als-variants}
\end{figure}

\Cref{fig:compare-nes-als-variants} shows the performance of our optimal Nesterov-ALS-RF-SG method 
on the synthetic test problems, with an ablation analysis that compares it with those
variants obtained by varying one hyper-parameter of Nesterov-ALS-RF-SG at a time.
In this $\tau$-plot, we display, for each method, the fraction of the 60 problem runs for which
the method execution time is within a factor $\tau$ of the fastest method for that problem.
For example, for $\tau=1$, the plot shows the fraction of the 60 problems for which each method
is the fastest. For $\tau=2$, the plot shows, for each method, the fraction of the 60 problems for which the method reaches the convergence tolerance in a time within a factor of two of the fastest method for that problem, etc.
As such, the area between curves is a measure for the relative performance of the methods, with the
curves at the top performing the best.

We can see that several variants have comparable performance to Nesterov-ALS-RF-SG, so
the optimal choice of restart mechanism and momentum weight is not very sensitive.
For these tests, changing the delay parameter has least effect on the performance.
Interestingly, this is then followed by changing the momentum weight to be a
constant of 1. It thus appears that, for our non-convex problem, the simple choice of $\beta_k=1$
combined with a suitable restart mechanism leads to a highly performant acceleration method.
This is followed by changing function restarting to gradient restarting and
speed restarting, respectively.
More detailed numerical results comparing Nesterov-ALS-RF-SG with a broader variation of
restarted Nesterov-ALS methods are shown in Appendix \ref{ap:detailed-results}, further confirming that
Nesterov-ALS-RF-SG generally performs the best among the family of restarted Nesterov-ALS
methods we have considered, for the synthetic test problems.
%\todo{line search momentum weight missing in this comparison. maybe include it}

\begin{figure}[h!]
%	\begin{subfigure}{\linewidth}
%		\includegraphics[width=\linewidth]{figs/New_Competetive_method_EasyALG1mod3.pdf}
%		\caption{Low collinearity}
%	\end{subfigure}
%	\begin{subfigure}{\linewidth}
	\centering{\stackunder{\includegraphics[width=.7\linewidth]{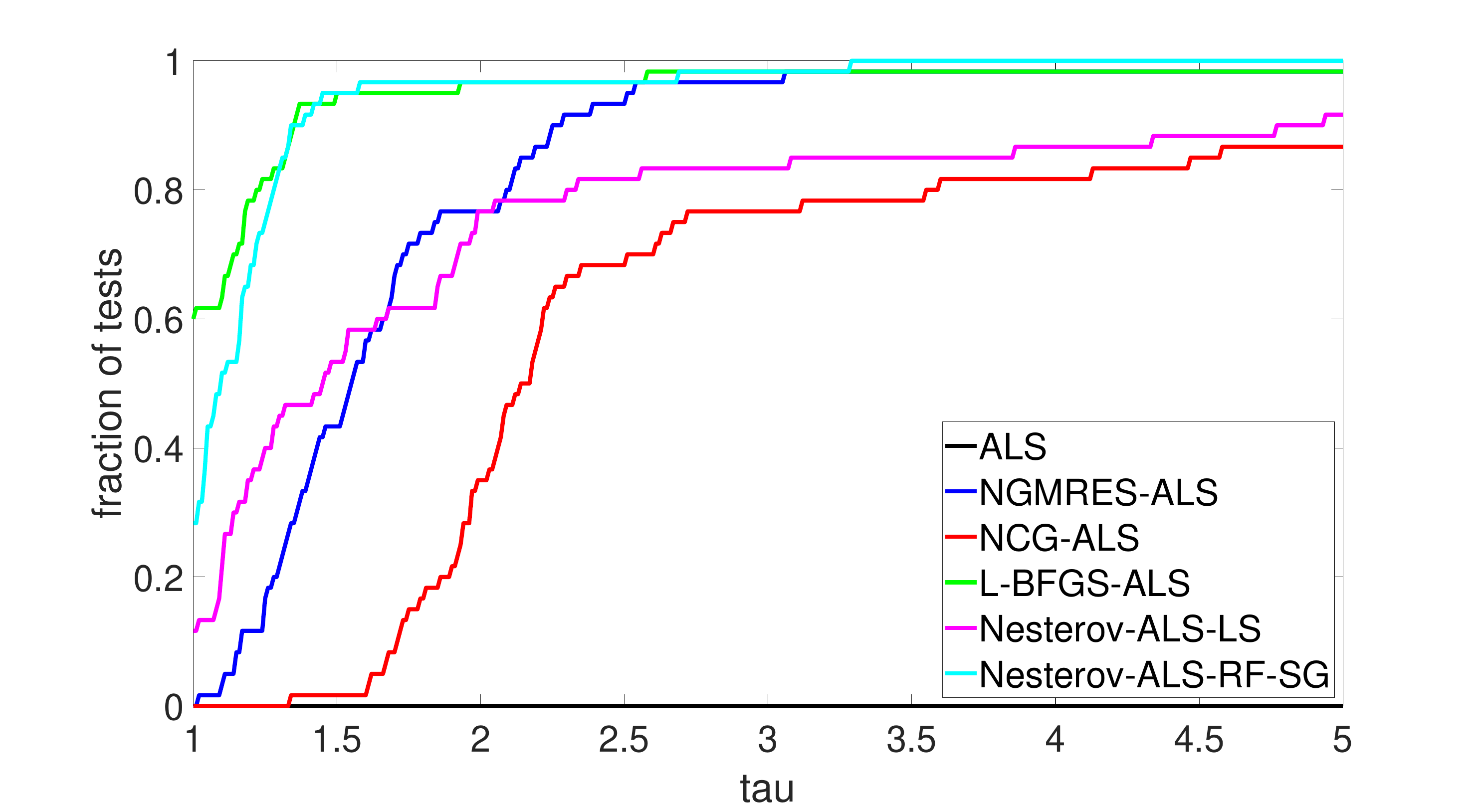}}
	{\rev{(a) $\tau$-plot for high accuracy tolerance value $tol=10^{-9}$.}}}
		\centering{\stackunder{\includegraphics[width=.7\linewidth]{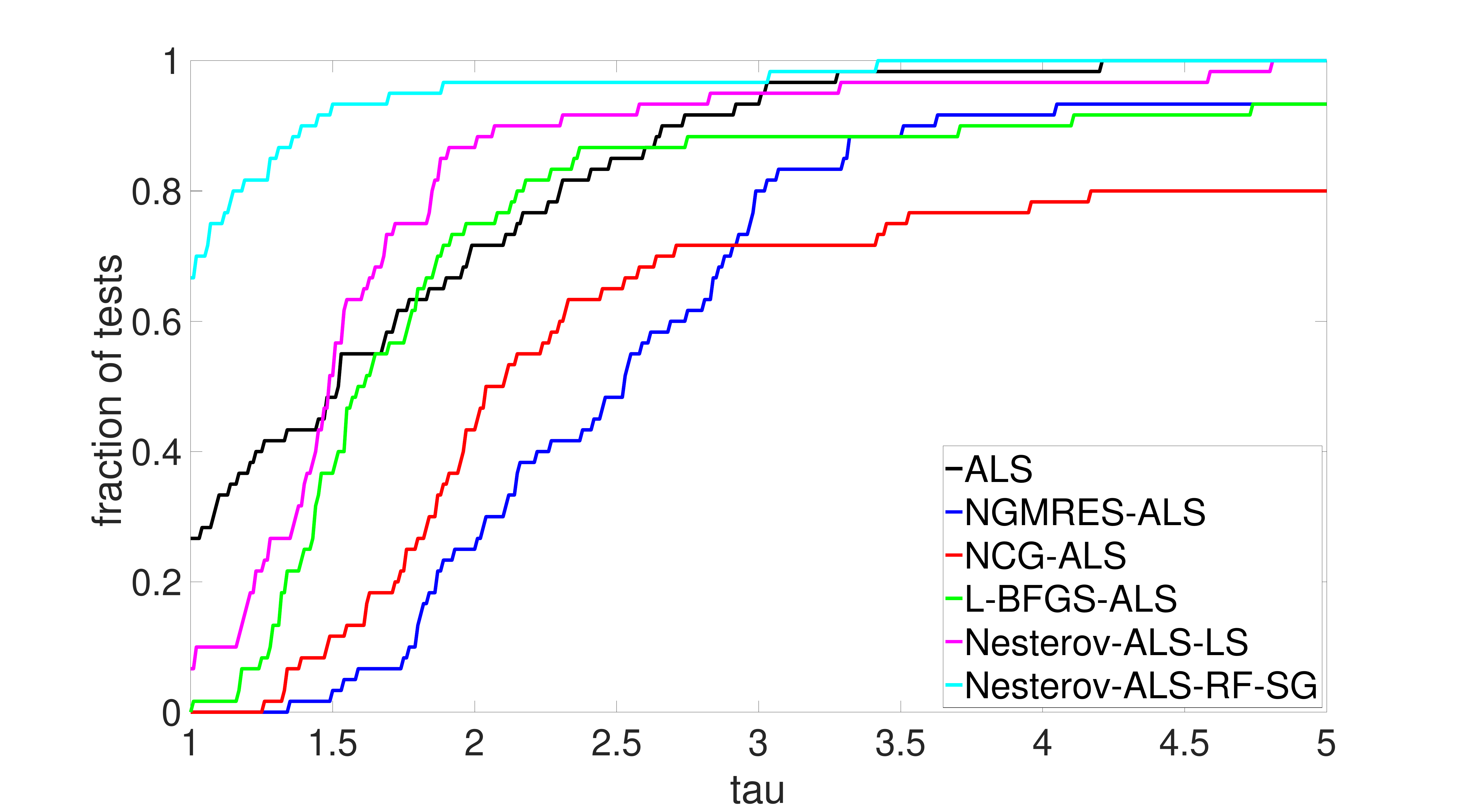}}
	{\rev{(b) $\tau$-plot for low accuracy tolerance value $tol=10^{-5}$.}}}
%		\caption{High collinearity}
%	\end{subfigure}
	\caption{\rev{Comparisons of the following algorithms on synthetic test problems:
	 the optimal restarted algorithm (Nesterov-ALS-RF-SG), the line search
	 restarted algorithm (Nesterov-ALS-LS), and existing accelerated ALS methods.}}
	\label{fig:compare-acc-methods-synthetic}
\end{figure}

\Cref{fig:compare-acc-methods-synthetic} compares Nesterov-ALS-RF-SG with the line
search version Nesterov-ALS-LS (equivalent or similar to existing line search methods to accelerate ALS\cite{harshman1970foundations,rajih2008enhanced,chen2011new,sorber2016exact}),
and with several other existing accelerated ALS methods, namely, GMRES-ALS, NCG-ALS, and LBFGS-ALS.
\rev{For high accuracy (top panel), Nesterov-ALS-RF-SG performs similarly to the
best performing of the existing methods we compare with, LBFGS-ALS\cite{sterck2018nonlinearly}. It performs substantially better than
Nesterov-ALS-LS (it avoids the expensive line searches). 
Nevertheless,
Nesterov-ALS-LS is competitive
with the existing NGMRES-ALS\cite{sterck2012nonlinear}, and superior to NCG-ALS\cite{sterck2015nonlinearly}.}

\rev{For low accuracy (bottom panel), Nesterov-ALS-RF-SG and Nesterov-ALS-LS still perform very well. ALS is now more competitive, which is not unexpected, since ALS is often efficient at reducing the initial error quickly, but then may converge slowly later on for difficult problems.}

%\medskip\noindent\emph{Effect of collinearity}.
%Collinearity was previously shown to have a strong influence on the performance
%of the ALS algorithm \cite{acar2011scalable}.

%Alternative_Methods_HardALG1mod3.pdf
%Broken_down_delay_onlyALG1mod3Hard.pdf
%Broken_down_restart_onlyALG1mod3Hard.pdf
%Broken_down_step_change_onlyALG1mod3Hard.pdf
%Paper draft.pdf
%char_plots_test_9_seed3_ALG1mod3Easy.pdf
%char_plots_test_9_seed9_ALG1mod3Easy.pdf
%char_plots_time_test_9_seed3_ALG1mod3Easy.pdf
%char_plots_time_test_9_seed9_ALG1mod3Easy.pdf
%plot2_HardALG1mod3.pdf

\begin{figure}[h!]
%	\begin{subfigure}{\linewidth}
	\centering{\stackunder{\includegraphics[width=.7\linewidth]{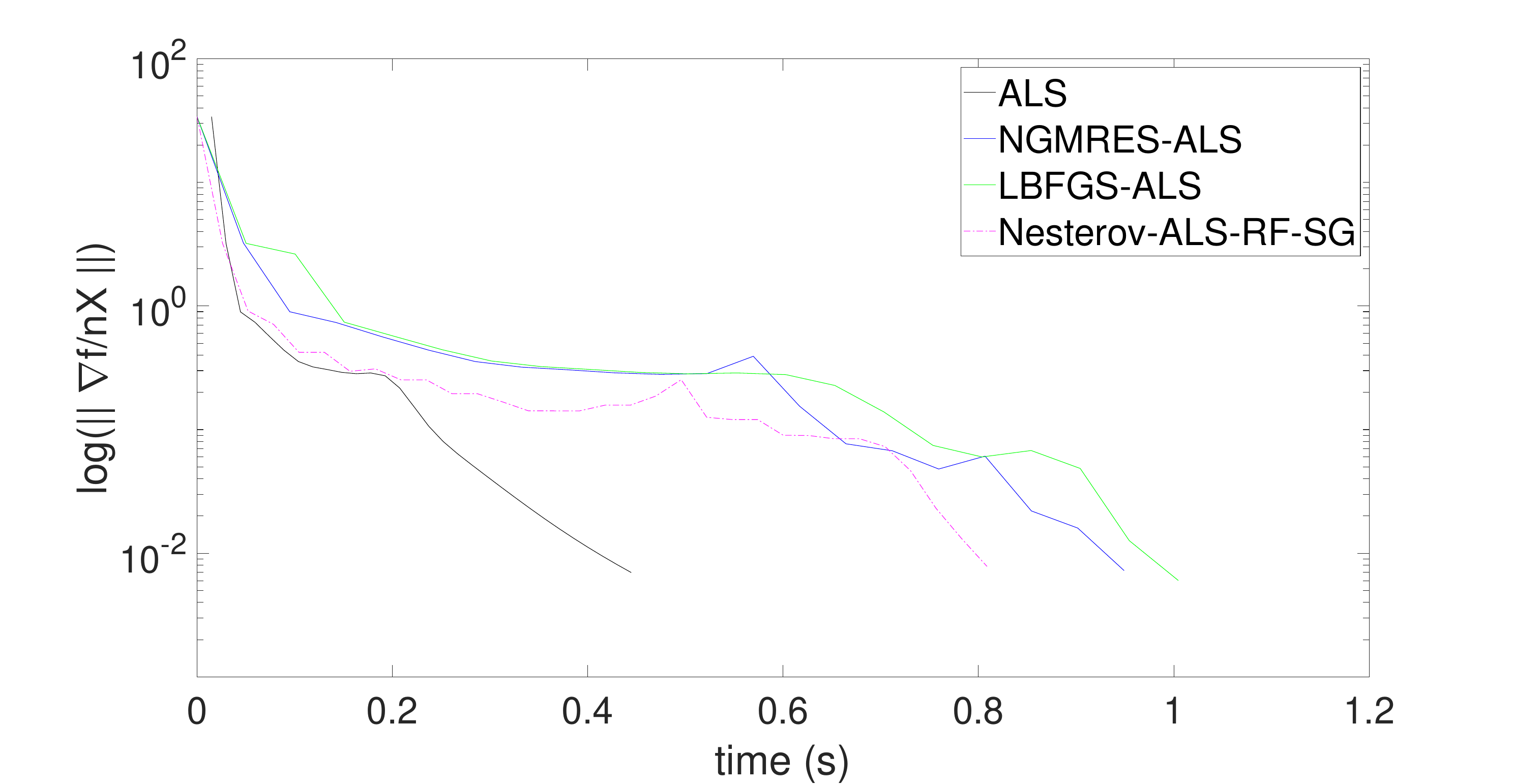}}
	{\rev{(a) Gradient convergence plot.}}}
%		\caption{time}
%	\end{subfigure}
%	\begin{subfigure}{\linewidth}
%		\includegraphics[width=\linewidth]{figs/New_enron_iter.pdf}
%		\caption{iteration}
%	\end{subfigure}
%	\begin{subfigure}{\linewidth}
		\centering{\stackunder{\includegraphics[width=.7\linewidth]{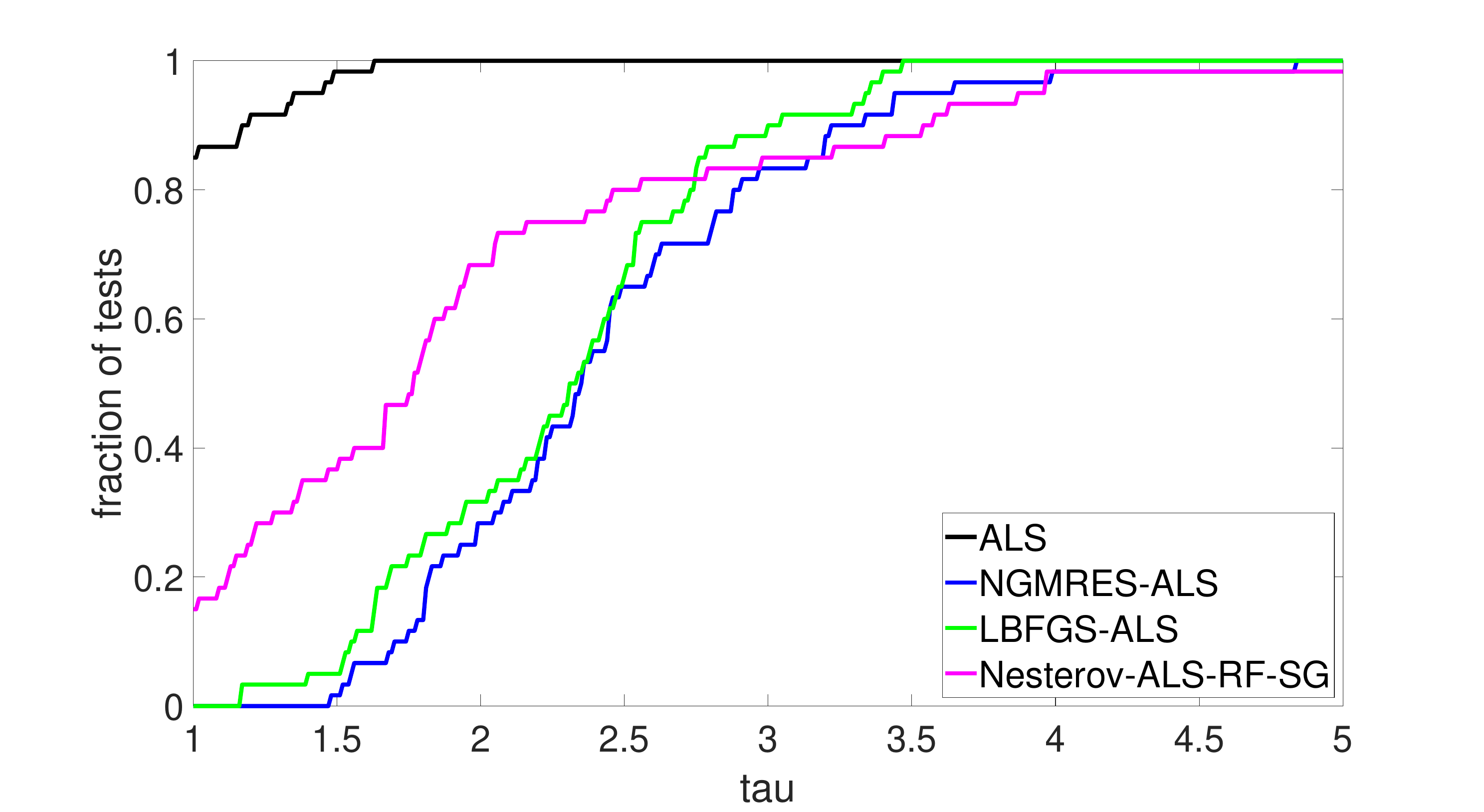}}
	{\rev{(b) $\tau$-plot for high accuracy tolerance value $tol=10^{-7}\,\|\nabla f(x_0)\|$.}}}
		\centering{\stackunder{\includegraphics[width=.7\linewidth]{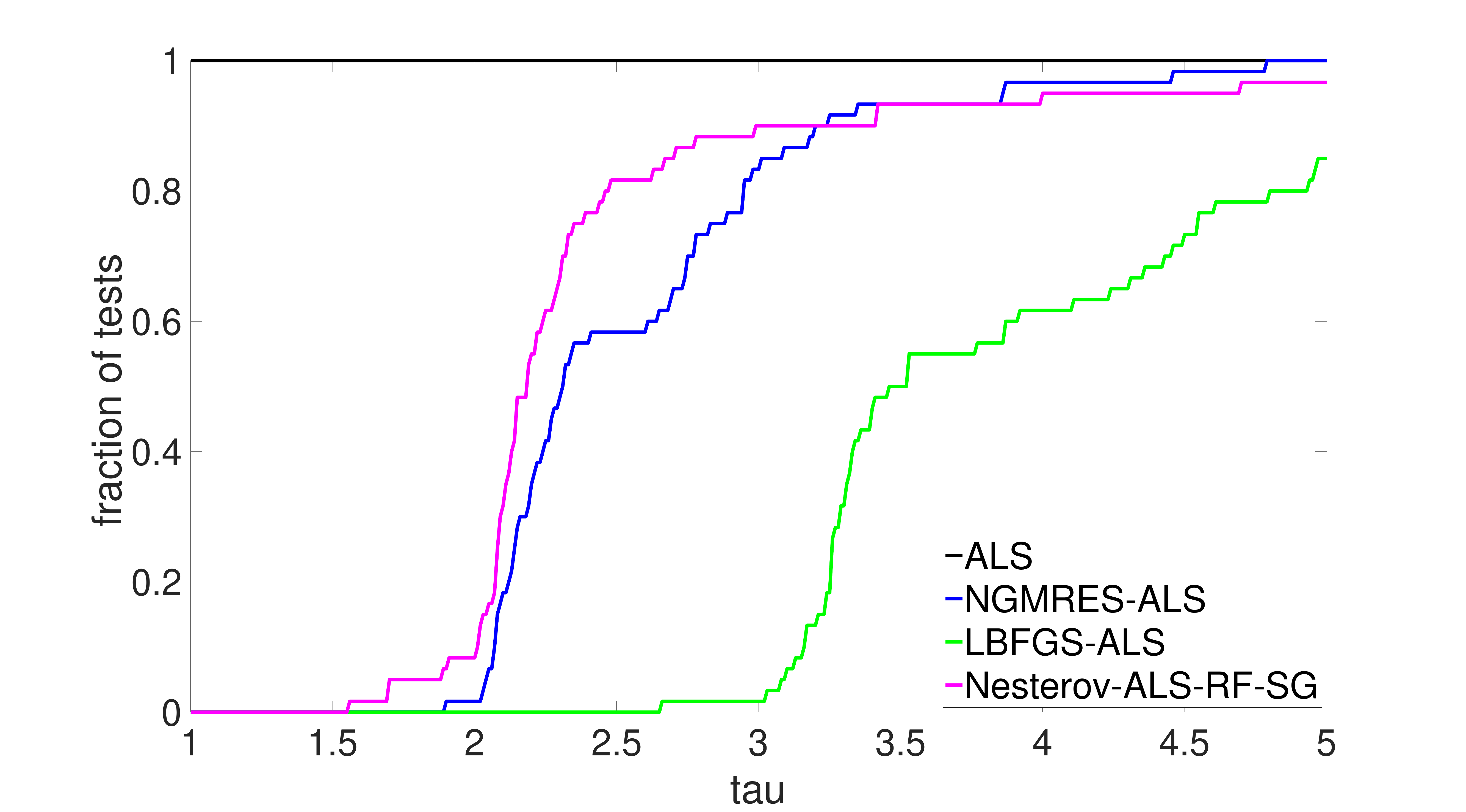}}
	{\rev{(c) $\tau$-plot for low accuracy tolerance value $tol=10^{-5}\,\|\nabla f(x_0)\|$.}}}
%		\caption{tau plot}
%	\end{subfigure}
	\caption{\rev{Comparison of algorithms on the Enron data.}}
	\label{fig:enron}
\end{figure}

\subsection{The Enron Dataset and Results.}
The Enron dataset is a subset of the corporate email communications that were
released to the public as part of the 2002 Federal Energy Regulatory Commission
(FERC) investigations following the Enron bankruptcy.
After various steps of pre-processing\cite{chi2012tensors},
a sender $\times$ receiver $\times$ month tensor of size 105$\times$105$\times$28 was obtained.
We perform rank-10 CP decompositions for Enron.
Fig.\ \ref{fig:enron} shows gradient norm convergence for one typical test run,
and a $\tau$-plot for 60 runs with different random initial guesses and 
convergence tolerances $tol=10^{-7}\,\|\nabla f(x_0)\|$ (high accuracy, middle
panel) and $tol=10^{-5}\,\|\nabla f(x_0)\|$ (lower accuracy, bottom panel).
For this well-conditioned problem, %(see discussion below), 
ALS converges fast and does not need acceleration. In fact, the acceleration overhead makes
ALS faster than any of the accelerated methods. This is consistent with previously reported results\cite{acar2011scalable,sterck2012nonlinear,sterck2015nonlinearly,sterck2018nonlinearly}
for well-conditioned problems.

\begin{figure}[h!]
%	\begin{subfigure}{\linewidth}
	\centering{\stackunder{\includegraphics[width=.7\linewidth]{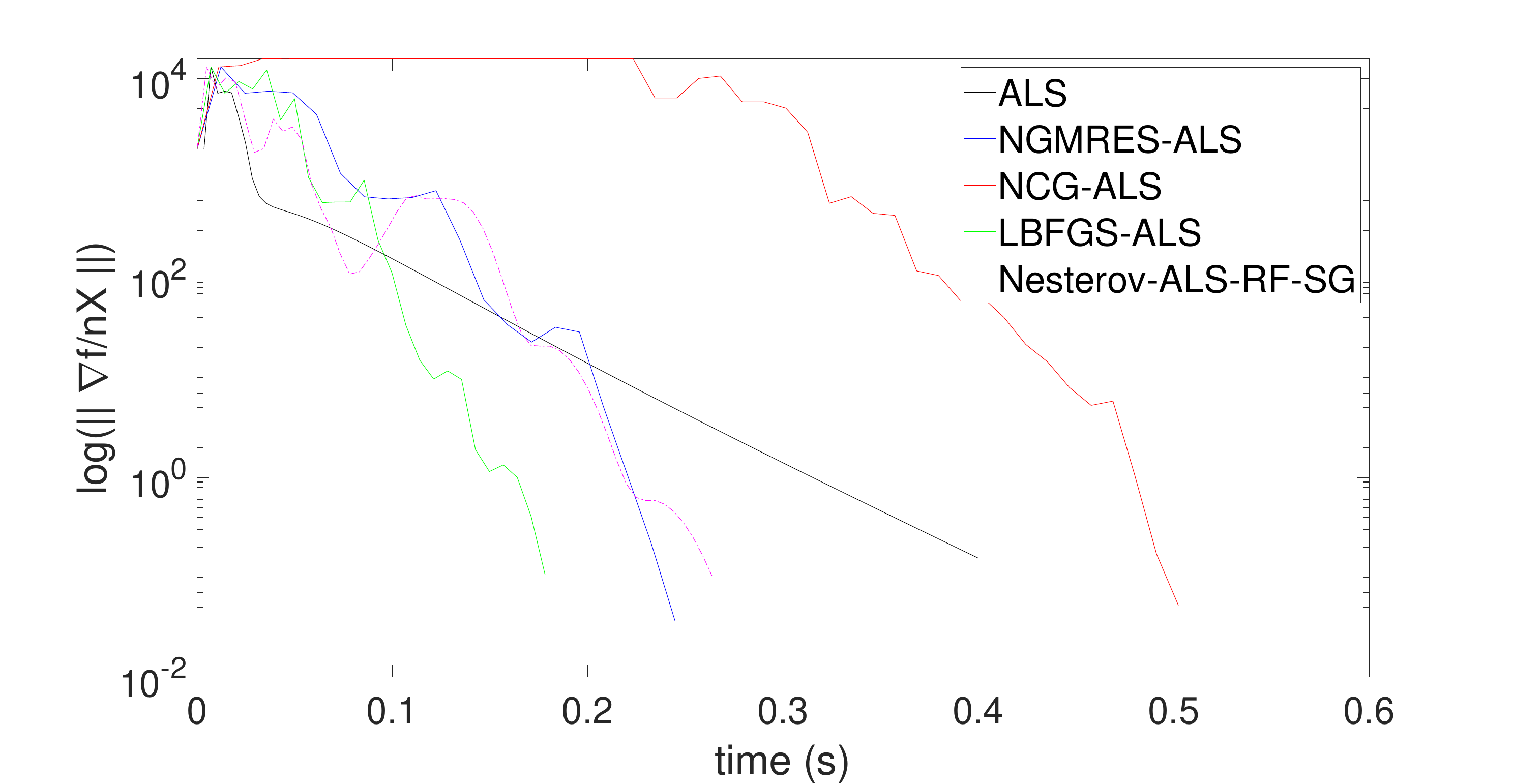}}	
	{\rev{(a) Gradient convergence plot.}}}
%		\caption{time}
%	\end{subfigure}
%	\begin{subfigure}{\linewidth}
%		\includegraphics[width=\linewidth]{figs/New_claus_iter.pdf}
%		\caption{iteration}
%	\end{subfigure}
%	\begin{subfigure}{\linewidth}
		\centering{\stackunder{\includegraphics[width=.7\linewidth]{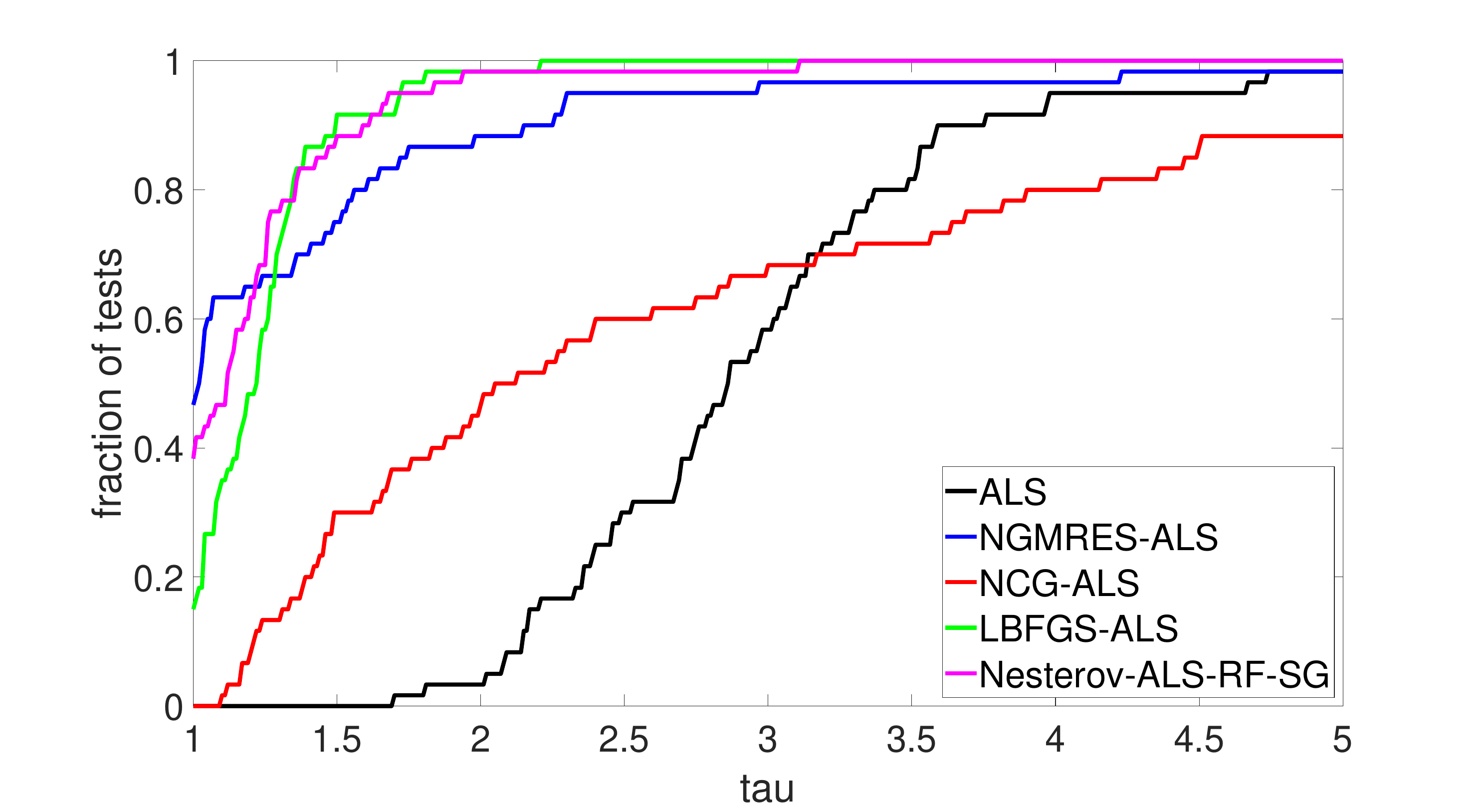}}
	{\rev{(b) $\tau$-plot for high accuracy tolerance value $tol=10^{-7}\,\|\nabla f(x_0)\|$.}}}
		\centering{\stackunder{\includegraphics[width=.7\linewidth]{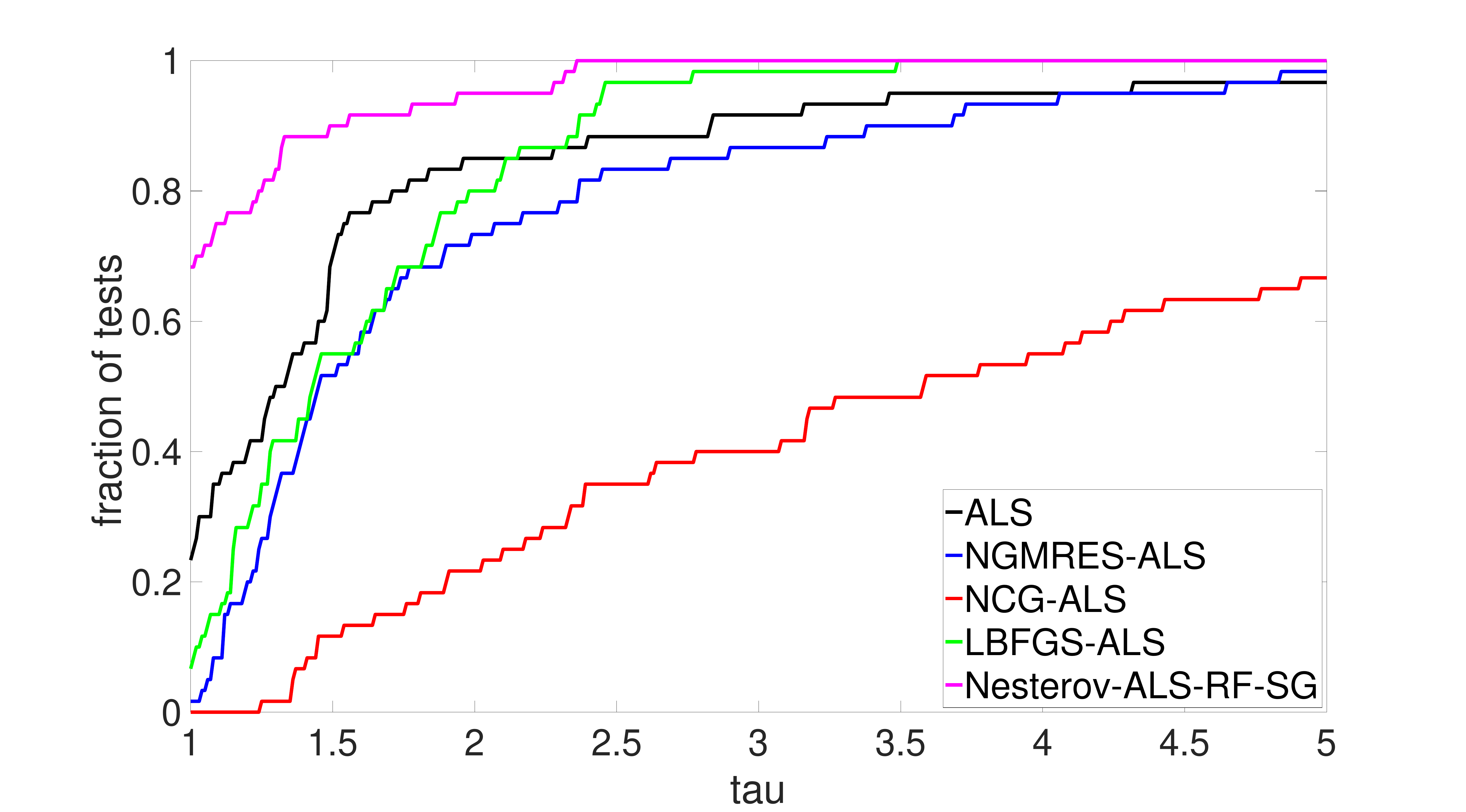}}
	{\rev{(c) $\tau$-plot for low accuracy tolerance value $tol=10^{-4}\,\|\nabla f(x_0)\|$.}}}
%		\caption{tau plot}
%	\end{subfigure}
	\caption{\rev{Comparison of algorithms on the Claus data.}}
	\label{fig:claus}
\end{figure}

\subsection{The Claus Dataset and Results.}
The claus dataset is a 5$\times$201$\times$61 tensor consisting of fluorescence measurements
of 5 samples containing 3 amino acids, taken for 201 emission wavelengths and
61 excitation wavelengths.
Each amino acid corresponds to a rank-one component\cite{andersen2003practical}. We perform a rank-3 CP decomposition for claus.
%Tryptophan (top), Tyrosine (middle), Phenylalanine (bottom)
Fig.\ \ref{fig:claus} shows gradient norm convergence for one test run,
and a $\tau$-plot for 60 runs with different random initial guesses and 
convergence tolerances $tol=10^{-7}\,\|\nabla f(x_0)\|$ (high accuracy, middle
panel) and $tol=10^{-4}\,\|\nabla f(x_0)\|$
(low accuracy, bottom panel).
For this medium-conditioned problem, substantial
acceleration of ALS can be obtained if high accuracy is desired,
and Nesterov-ALS-RF-SG performs
as well as the best methods we compare with, but it is much easier to implement.
\rev{For low accuracy, ALS is more competitive, but Nesterov-ALS-RF-SG still outperforms it.}

\begin{figure}[h!]
%	\begin{subfigure}{\linewidth}
	\centering{\stackunder{\includegraphics[width=.7\linewidth]{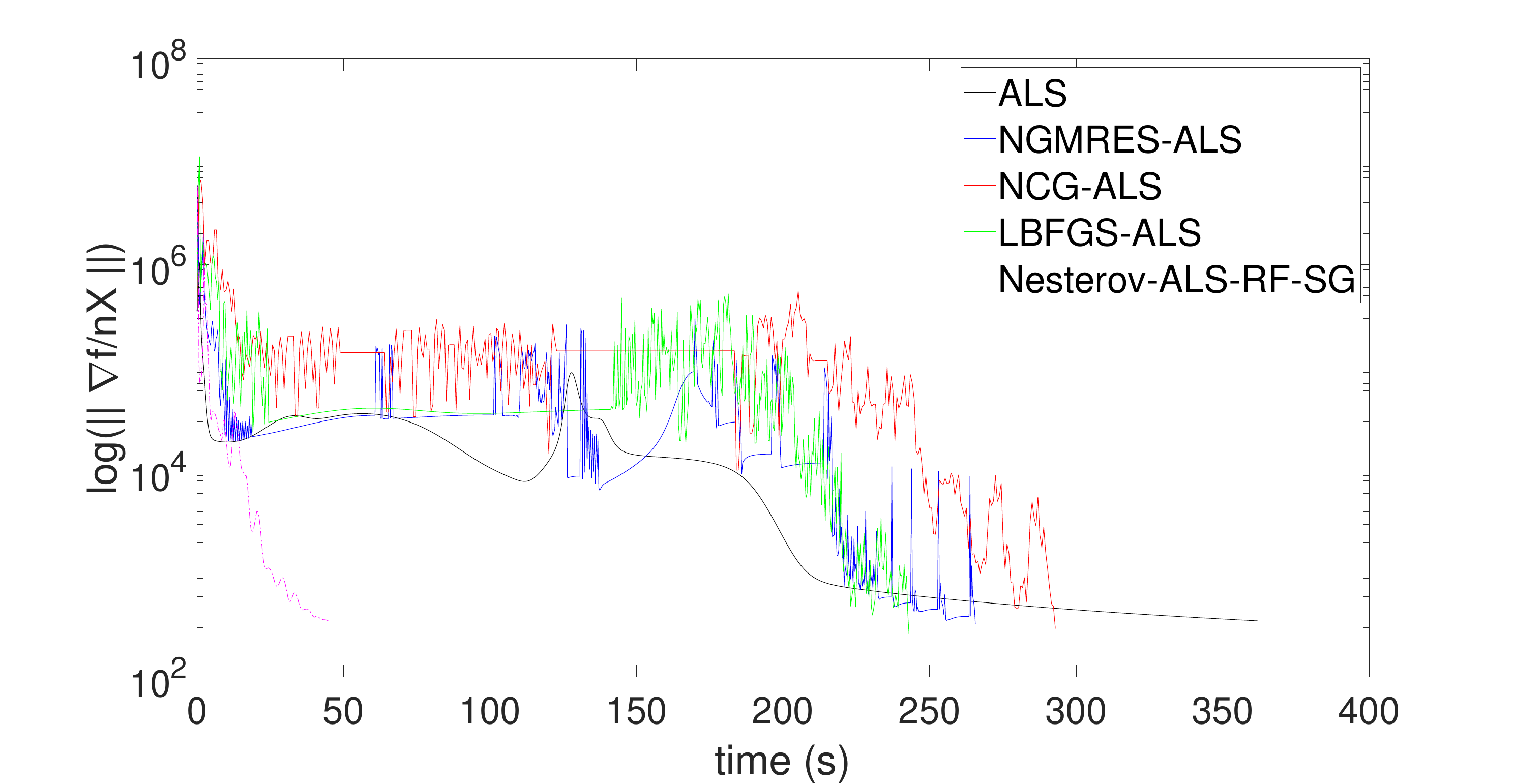}}
	{\rev{(a) Gradient convergence plot.}}}
%		\caption{time}
%	\end{subfigure}
%	\begin{subfigure}{\linewidth}
%		\includegraphics[width=\linewidth]{figs/New_gas3_iter.pdf}
%		\caption{iteration}
%	\end{subfigure}
%	\begin{subfigure}{\linewidth}
		\centering{\stackunder{\includegraphics[width=.7\linewidth]{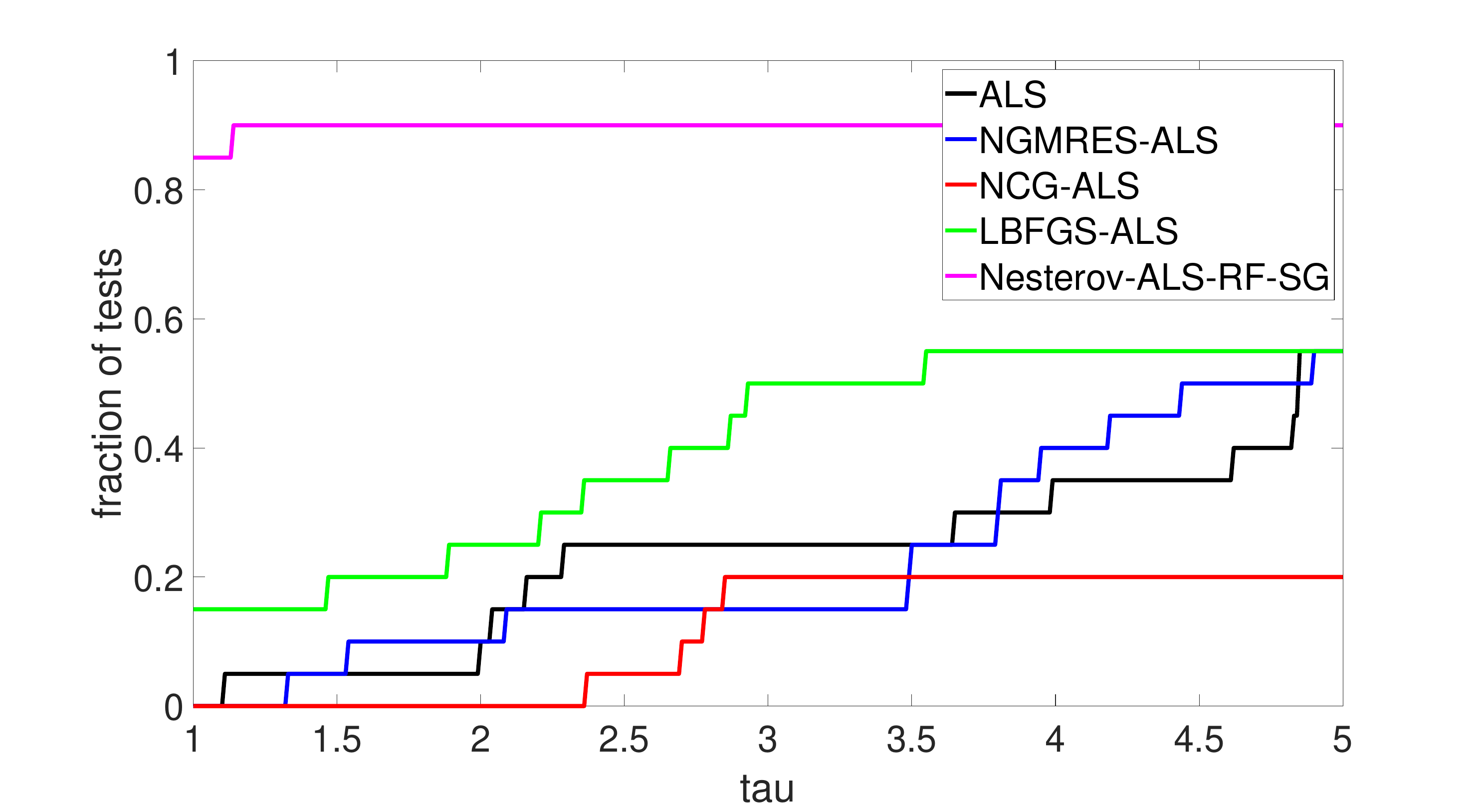}}
		{\rev{(b) $\tau$-plot for high accuracy tolerance value $tol=10^{-7}\,\|\nabla f(x_0)\|$.}}}
		\centering{\stackunder{\includegraphics[width=.7\linewidth]{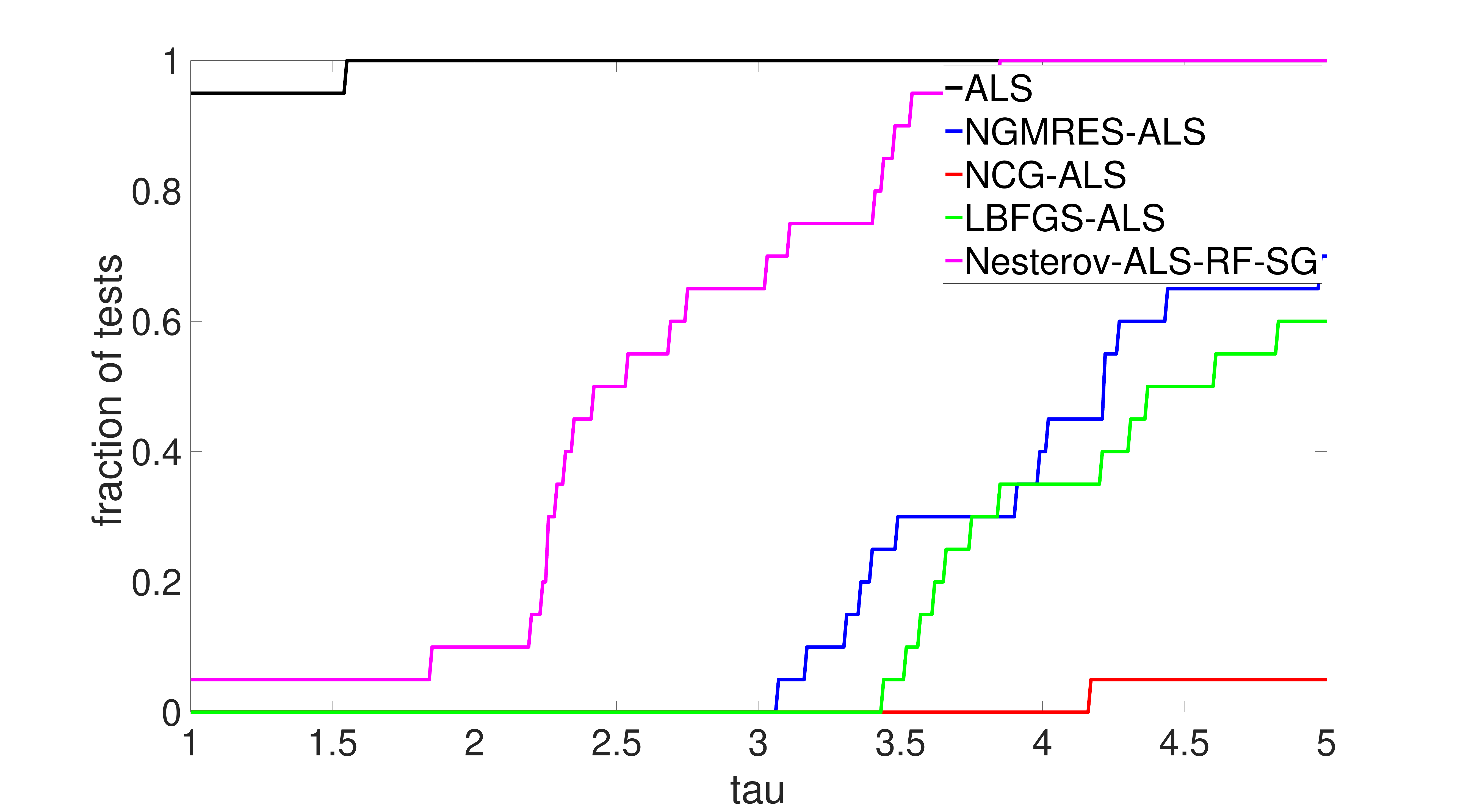}}
		{\rev{(c) $\tau$-plot for low accuracy tolerance value $tol=10^{-4}\,\|\nabla f(x_0)\|$.}}}
		
%		\caption{tau plot}
%	\end{subfigure}
	\caption{\rev{Comparison of algorithms on the Gas3 data.}}
	\label{fig:gas3}
\end{figure}

\subsection{The Gas3 Dataset and Results.}
Gas3 is relatively large and has multiway structure. 
It is a 71$\times$1000$\times$900 tensor consisting of readings from 71 chemical sensors used
for tracking hazardous gases over 1000 time steps\cite{vergara2013performance}.
%In each experiment, a gas was released into a chamber, and the readings of 71
%sensors were recorded for 1000 time steps. 
There were three gases, and 300 experiments were performed for each gas, varying fan
speed and room temperature. We perform a rank-5 CP decomposition for Gas3.

Fig.\ \ref{fig:gas3} shows gradient norm convergence for one typical test run,
and a $\tau$-plot for 20 runs with different random initial guesses and 
convergence tolerance $tol=10^{-7}\,\|\nabla f(x_0)\|$ (high accuracy, middle panel) and $tol=10^{-4}\,\|\nabla f(x_0)$
(low accuracy, bottom panel).

For this highly ill-conditioned problem, 
ALS converges slowly, and NGMRES-ALS, NCG-ALS and LBFGS-ALS
behave erratically. %NCG-ALS is not included in the $\tau$-plot because it does not converge.
For high accuracy, our newly proposed Nesterov-ALS-RF-SG very substantially
outperforms all other methods (not only for the convergence profile shown in the
top panel, but for the large majority of the 20 tests with random 
initial guesses). Nesterov-ALS-RF-SG performs much more
robustly for this highly ill-conditioned problem than
any of the other accelerated methods, and reaches high accuracy
much faster than any other method.
We were initially surprised that Nesterov-ALS-RF-SG performs so much better than the
other accelerated ALS methods we compare with, and have so far not found a clear
indication why this is the case. One possible explanation is that the line search employed in NGMRES-ALS,
NCG-ALS and LBFGS-ALS may suffer robustness issues due to the ill-conditioning of the problem.
\rev{For low accuracy, ALS is clearly the fastest: it is efficient at reducing the initial error quickly, but converges slowly later on for this difficult problem.}

The other tests in this paper indicate that our proposed Nesterov-ALS-RF-SG acceleration method
is competitive with leading existing acceleration methods for ALS, and, additionally, this result gives some initial indication
that Nesterov-ALS-RF-SG is surprisingly robust for highly ill-conditioned problems. This will need to be investigated
further when these methods, which will be made available as part of the Poblano toolbox for optimization,
will be applied by us and others to other demanding CP decomposition problems.

%\subsection{Discussion on Real-World Problems.}
%We speculated that our accelerated ALS methods may work best for ill-conditioned problems.
%To verify this, we computed the condition number of the initial Hessians for the
%three real-world problems. These were 58,842, 3,094,000, and 119,220,000 for Enron, Claus, and Gas3,
%respectively. This agrees with the observed advantage of Nesterov-ALS-RF-SG in 
%\Cref{fig:enron,fig:claus,fig:gas3}.

%Nesterov-ALS-LS requires fewest iterations in general? 
%If yes, this suggests that developed customized line search technique (like
%the one mentioned in an email) is worthwhile.
%Nesterov-ALS with momentum weight 1?

%---------------------------------------------------------------------
\subsection{\rev{Discussion on Per-Iteration Computational Cost}}
%---------------------------------------------------------------------
\rev{It is important to compare the computational cost per iteration of our
Nesterov-accelerated ALS methods with plain ALS iterations. For simplicity, we
discuss the case of dense order-3 tensors of size $I \times J \times K$.
For $R \ll I,J,K$, the dominant cost in each ALS step is the formation of 
the three MTTKRPs on the left-hand side in Eqs.\
(\ref{eq:ALS1})-(\ref{eq:ALS3}), which takes approximately $O(IJKR)$
arithmetic operations.
%Thus the time complexity of an ALS update is $O(IJKR)$. 
At the end of an ALS iteration, the function value can be computed at virtually
no cost by noting that 
\begin{align*}
	f(A,B,C)
	= \frac{1}{2}\norm{T - \widetilde{T}}_{F}^2
	= \frac{1}{2}\norm{T_3-C (A \odot B)^T}_F^2
	= \frac{1}{2}\norm{T}_{F}^{2} - T_{3} (A \odot B) \cdot C + \frac{1}{2} \norm{C (A \odot B)^{T}}_{F}^{2},
\end{align*}
where $\norm{T}_F^2$ is pre-computed, the MTTKRP $T_{3} (A \odot B)$ is re-used,
and $\cdot$ is the dot product between two matrices (i.e., the sum of the entries of the element-wise multiplication of the matrices). The third term equals $\frac{1}{2}\|\widetilde{T}\|_F^2$, and, due to the structure of the rank-1 terms, it can be computed in $O(2R^2(I+J+K))$ operations, which is negligible compared to the cost of the MTTKRPs in ALS.
%\cite{phan2012fast} 
In fact, this efficient computation is the default algorithm for computing the function value in the ALS implementation of the Tensor Toolbox\cite{TTB_Software}.}

\rev{Computing the three MTTKRPs is also the dominant cost when calculating the gradient of $f$, according to Eqs.\ (\ref{eq:grad1})-(\ref{eq:grad3}), and following a gradient computation the function value comes virtually for free, so the dominant cost of a function + gradient evaluation is also of $O(IJKR)$, the same cost as an ALS iteration.}

\rev{In our Nesterov-accelerated ALS methods, we use the gradient and/or function value to compute the restart condition and/or the step length, so we perform one function + gradient evaluation per iteration in the implementation of all our methods, in addition to the cost of one ALS step per iteration. For this reason, each of our Nesterov-accelerated ALS iterations is about twice as expensive as a regular ALS iteration. Still, the numerical results in the previous sections show clearly that the Nesterov acceleration tends to reduce the number of iterations required for (accurate) convergence by much more than a factor of 2, making the accelerated methods clear winners over ALS (despite the doubled cost per iteration), for difficult problems or when accurate solutions are required.}

\rev{Note that it is possible to avoid the doubling of the per-iteration cost by considering acceleration variants that only rely on function values for the restart mechanism, and determine the step length without using gradient information. In that respect, our accelerated method with function restart and constant momentum weight one is attractive. Note also that the restarted acceleration methods in Ang and Gillis\cite{ang2019accelerating} for Nonnegative Matrix Factorization use more elaborate strategies for determining step lengths that also do not require gradient information.}

\rev{Finally, it is interesting to observe that, due to the line searches, the per-iteration cost of our NGMRES-ALS, Anderson-ALS and LBFGS-ALS methods is typically 4 to 5 times the cost of 1 ALS iteration, and the per-iteration cost of NCG-ALS is typically 6 to 7 times the cost of 1 ALS iteration. Still, some of these methods achieve very large reductions in the number of iterations required for convergence, and are often among the most efficient, in particular, LBFGS-ALS and NGMRES-ALS/Anderson-ALS.}

\begin{figure}[h!]
	\centering{\stackunder{\includegraphics[width=.6\linewidth]{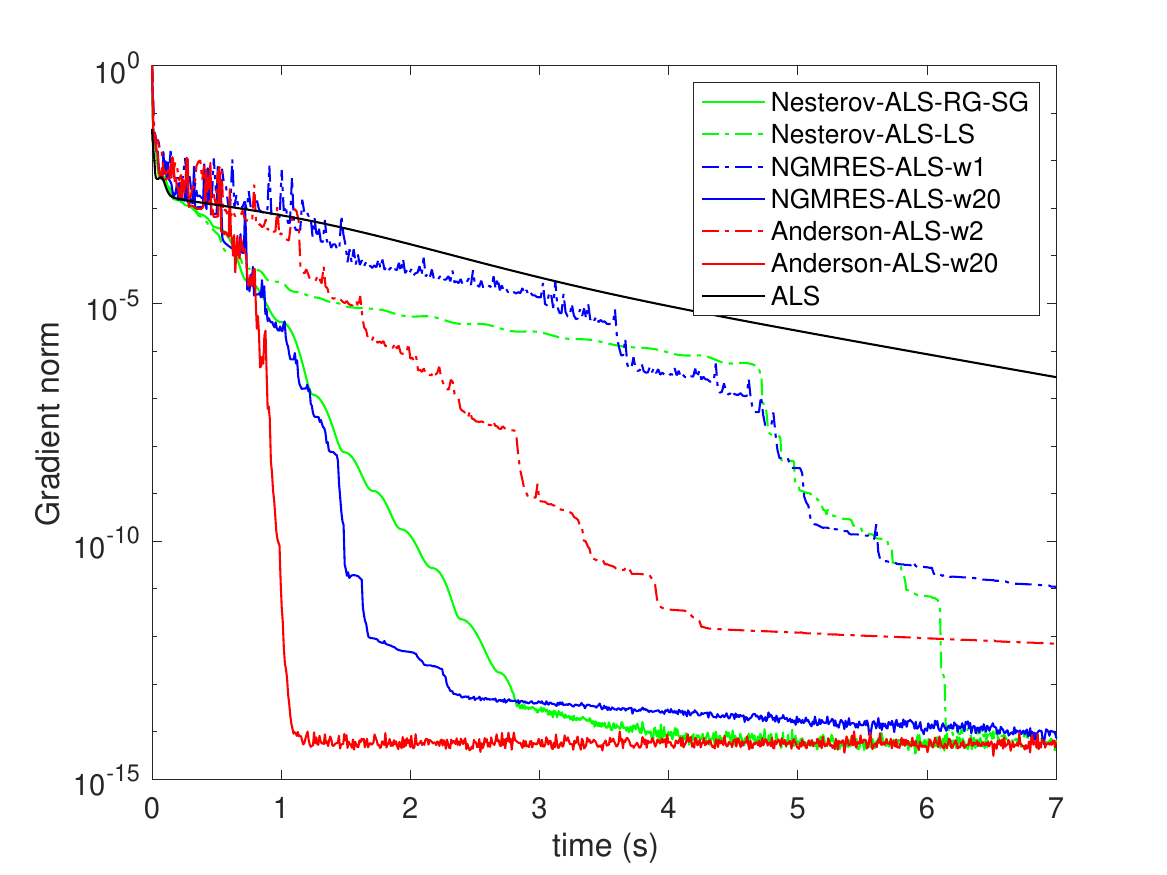}}
	{(a) Gradient convergence plot.}}
		\centering{\stackunder{\includegraphics[width=.6\linewidth]{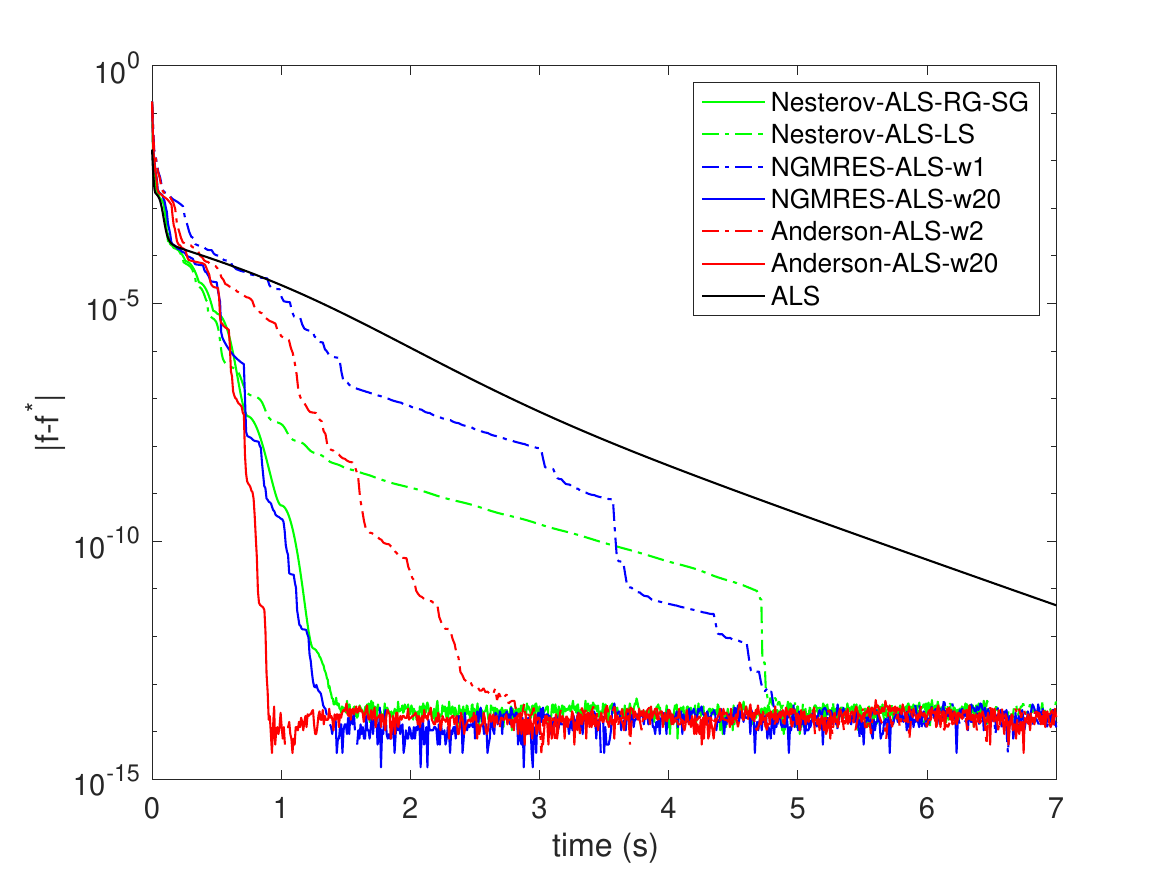}}
	{(b) Function value convergence plot.}}
	\caption{\rev{Convergence comparison with existing line search methods for synthetic test problem 2 from Table \ref{tab:synthetic}. Top panel: gradient norm as a function of time. Bottom panel: convergence of the objective, $f$, to its minimum value, $f^*$, as a function of time.
One of our proposed Nesterov-accelerated methods with restart, Nesterov-ALS-RG-SG (green), is compared with Harshman's extrapolation with line search\cite{harshman1970foundations} (equivalent to NGMRES-ALS-w1, blue-dashed), Chen et al.'s 1-step extrapolation with line search\cite{chen2011new} (equivalent to Nesterov-ALS-LS, green-dashed), and Chen et al.'s 2-step extrapolation with line search\cite{chen2011new} (similar to our Anderson acceleration with window size 2, Anderson-ALS-w2, red-dashed). Our proposed restarted Nesterov method (green) converges substantially faster than the existing line search methods (dashed). Also, for Harshman's extrapolation with line search (blue-dashed), the result for NGMRES-ALS-w20 (blue) shows that a window size for NGMRES greater than 2 leads to much faster convergence. Similarly, for Chen et al.'s 2-step extrapolation with line search (similar to Anderson-ALS-w2, red-dashed, the result for Anderson-ALS-w20 (red) shows that a window size for Anderson acceleration greater than 2 leads to much faster convergence. These convergence plots are for a specific example, but they are typical for the relative performance of the methods.}}
	\label{fig:compareLS}
\end{figure}

%---------------------------------------------------------------------
\subsection{\rev{Further Comparison with Line Search Methods}}
%---------------------------------------------------------------------

\rev{
It is interesting to further compare numerical results with the existing line search methods that were discussed in Section 
\ref{subsec:compare-line-search}.
Figure \ref{fig:compareLS} compares convergence for synthetic test problem 2 from Table \ref{tab:synthetic}, for several methods. 
The top panel shows the gradient norm as a function of time, and the bottom panel the convergence of the objective, $f$, to its minimum value, $f^*$, as a function of time. We compute 
$f^*$ as the smallest value of $f$ obtained by any algorithm after performing sufficiently many iterations to make $\nabla f$ converge to machine accuracy.}

\rev{One of our proposed Nesterov-accelerated methods with restart, Nesterov-ALS-RG-SG (green), is compared with several existing line search methods.}

\rev{We first compare with NGMRES-ALS with window size 1, NGMRES-ALS-w1 (blue-dashed), which is equivalent to Harshman's extrapolation with line search\cite{harshman1970foundations,rajih2008enhanced}. We use the Mor\'{e}-Thuente cubic line search rather than an exact line search\cite{rajih2008enhanced}, but it is known that the exact line search is expensive since it requires $2N - 1$ function evaluations\cite{acar2011scalable,sterck2012nonlinear}, whereas the cubic line search typically only requires about 2 or 3 function evaluations\cite{acar2011scalable}.}

\rev{We also compare with Nesterov extrapolation with line search, Nesterov-ALS-LS (green-dashed), which is equivalent with
Chen et al.'s 1-step extrapolation with line search\cite{chen2011new}, except that we don't use an exact line search.
Finally, we compare with Anderson acceleration with window size 2, Anderson-ALS-w2 (red-dashed), which uses the same
form of extrapolation as Chen et al.'s 2-step extrapolation\cite{chen2011new}, but determines optimal expansion coefficients by solving a least-squares problem in each step, whereas Chen et al.'s method makes an ad-hoc choice for the expansion coefficients that is the same in all iterations. As such, we can use Anderson-ALS-w2 as an (optimistic) proxy for estimating the performance of Chen et al.'s 2-step extrapolation\cite{chen2011new}. }

\rev{The results in Figure \ref{fig:compareLS} show that our proposed Nesterov-accelerated method with restart, Nesterov-ALS-RG-SG (green), converges much faster than the three existing line search methods (dashed). Also, for Harshman's extrapolation with line search (blue-dashed, NGMRES-ALS-w1), the result for NGMRES-ALS-w20 (blue) shows that a window size for NGMRES greater than 2 leads to much faster convergence. Similarly, for Chen et al.'s 2-step extrapolation with line search (similar to Anderson-ALS-w2, red-dashed), the result for Anderson-ALS-w20 (red) shows that a window size for Anderson acceleration greater than 2 leads to much faster convergence. These convergence plots are for a specific example, but they are typical for the relative performance of the methods.}

\rev{Note, finally, that the problem in Figure \ref{fig:compareLS} is an ill-conditioned problem, so $|f-f^*|$ reaches machine accuracy before the gradient does.}

%%%%%%%%%%%%%%%%%%%%%%%%%%%%%%%%%%%%%
\section{Conclusion} \label{sec:conclude}
%%%%%%%%%%%%%%%%%%%%%%%%%%%%%%%%%%%%%
We have proposed Nesterov-ALS methods with effective choices for momentum weights and restart mechanisms
that are simple and easy to implement as compared to several existing nonlinearly accelerated ALS methods,
such as GMRES-ALS, NCG-ALS, and LBFGS-ALS\cite{sterck2012nonlinear,sterck2015nonlinearly,sterck2018nonlinearly}. %It incurs little computational overhead on top of the vanilla ALS.
The optimal variant, using function restarting and gradient ratio momentum
weight, is competitive with or superior to stand-alone ALS and GMRES-ALS, NCG-ALS, LBFGS-ALS and line search ALS\cite{harshman1970foundations,rajih2008enhanced,chen2011new,sorber2016exact}.

Simple nonlinear iterative optimization methods like ALS and coordinate descent (CD)
are widely used in a variety of application domains.
There is clear potential for extending our approach to accelerating such simple optimization
methods for other non-convex problems. A specific example is Tucker tensor decomposition\cite{kolda2009tensor}. 
NCG, NGMRES and LBFGS acceleration have been applied to Tucker decomposition\cite{Hans_Tucker_decomp,sterck2018nonlinearly}, using a manifold approach to maintain the Tucker orthogonality constraints, and this approach can directly be extended to the Nesterov acceleration methods discussed in this paper.

More generally, we have formulated a Nesterov-type acceleration approach that can effectively accelerate optimization algorithms different from gradient descent (such as ALS) and for non-convex problems (such as CP tensor decomposition), using simple choices for momentum weights (as simple as setting them equal to one) and suitable restart mechanisms. The proposed methods are competitive in terms of performance with any existing acceleration methods for ALS and are very simple to implement, and initial tests on a challenging real-world problem indicate desirable robustness properties.

\rev{Matlab code implementing the proposed Nesterov acceleration methods is freely available at
\url{https://github.com/hansdesterck/nonlinear-preconditioning-for-optimization}.
The code includes implementations of all Nesterov-ALS variants, as well as NCG-ALS, NMGRES-ALS, LBFGS-ALS and Anderson-ALS.
The implementations are generic in that any suitable nonlinear preconditioner can be provided, not just ALS, and the nonlinearly preconditioned methods can be applied to any suitable optimization problem, not just canonical tensor decomposition.}

%
%%%%%%%%%%%%%%%%%%%%%%%%%%%%%%%%
%

%\acks
\section*{Acknowledgments}
\rev{The work of H. D. S. was partially supported by an NSERC Discovery Grant.}

\bibliographystyle{ieeetr}
\bibliography{ref}

%%%%%%%%%%%%%%%%%%%%%%%%%%%%%%%%%%%%%%%%
%%%%%%%%%%%%%%%%%%%%%%%%%%%%%%%%%%%%%%%%
%\newpage
%\clearpage
%
%%%%%%%%%%%%%%%%%%%%%%%%%%%%%%%%
%
\appendix

\section{Parameters for synthetic CP test problems}\label{ap:params}

Table \ref{table:test_problems} lists the parameters for the standard ill-conditioned synthetic test problems used in the paper\cite{acar2011scalable}.
The specific choices of parameters
for the six classes in Table \ref{table:test_problems} correspond to test problems 7-12 in De Sterck\cite{sterck2012nonlinear}. All tensors have equal size
$s=I_1=I_2=I_3$ in the three tensor dimensions, and have high collinearity $c$.
The six classes differ in their choice of tensor sizes ($s$), decomposition rank
($R$), and noise parameters $l_1$ and $l_2$.

\begin{table}[h!]
	\centering
	\caption{List of parameters for synthetic CP test problems.\label{table:test_problems}
}
	\begin{tabular}{|c | c | c | c | c | c | }
		\hline
		problem & $s$ & $c$ & $R$ & $l_1$ & $l_2$ \\
		\hline\hline
%		1 & 20 & 0.5 & 3 &  1 & 1\\
%		\hline
%		2 & 20 & 0.5 & 5 &  10 & 5\\
%		\hline
%		3 & 50 & 0.5 & 3 &  1 & 1\\
%		\hline
%		4 & 50 & 0.5 & 5 &  10 & 5\\
%		\hline
%		5 & 100 & 0.5 & 3 &  1 & 1\\
%		\hline
%		6 & 100 & 0.5 & 5 &  10 & 5\\
%		\hline\hline
		1 & 20 & 0.9 & 3 &  0 & 0\\
		\hline
		2 & 20 & 0.9 & 5 &  1 & 1\\
		\hline
		3 & 50 & 0.9 & 3 &  0 & 0\\
		\hline
		4 & 50 & 0.9 & 5 &  1 & 1\\
		\hline
		5 & 100 & 0.9 & 3 &  0 & 0\\
		\hline
		6 & 100 & 0.9 & 5 &  1 & 1\\
		\hline
	\end{tabular}\label{tab:synthetic}
\end{table}

\medskip

\section{Detailed comparisons for different restarting strategies}\label{ap:detailed-results}

\begin{figure*}[tb]
	\begin{centering}
	\begin{subfigure}{.8\linewidth}
		\centering{\includegraphics[width=.87\linewidth]{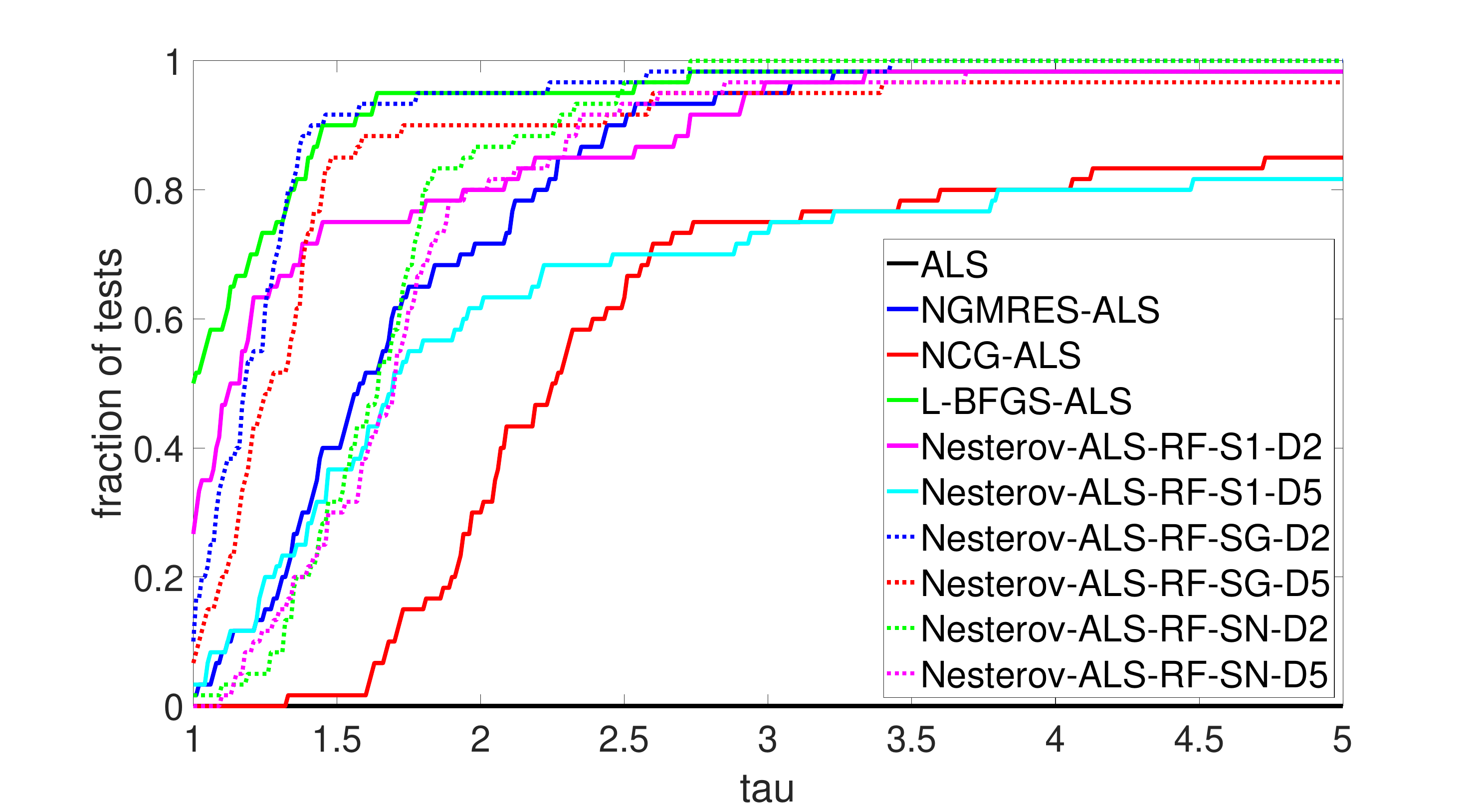}}
		\caption{variants with delay}
	\end{subfigure}
	\begin{subfigure}{.8\linewidth}
		\centering{\includegraphics[width=.87\linewidth]{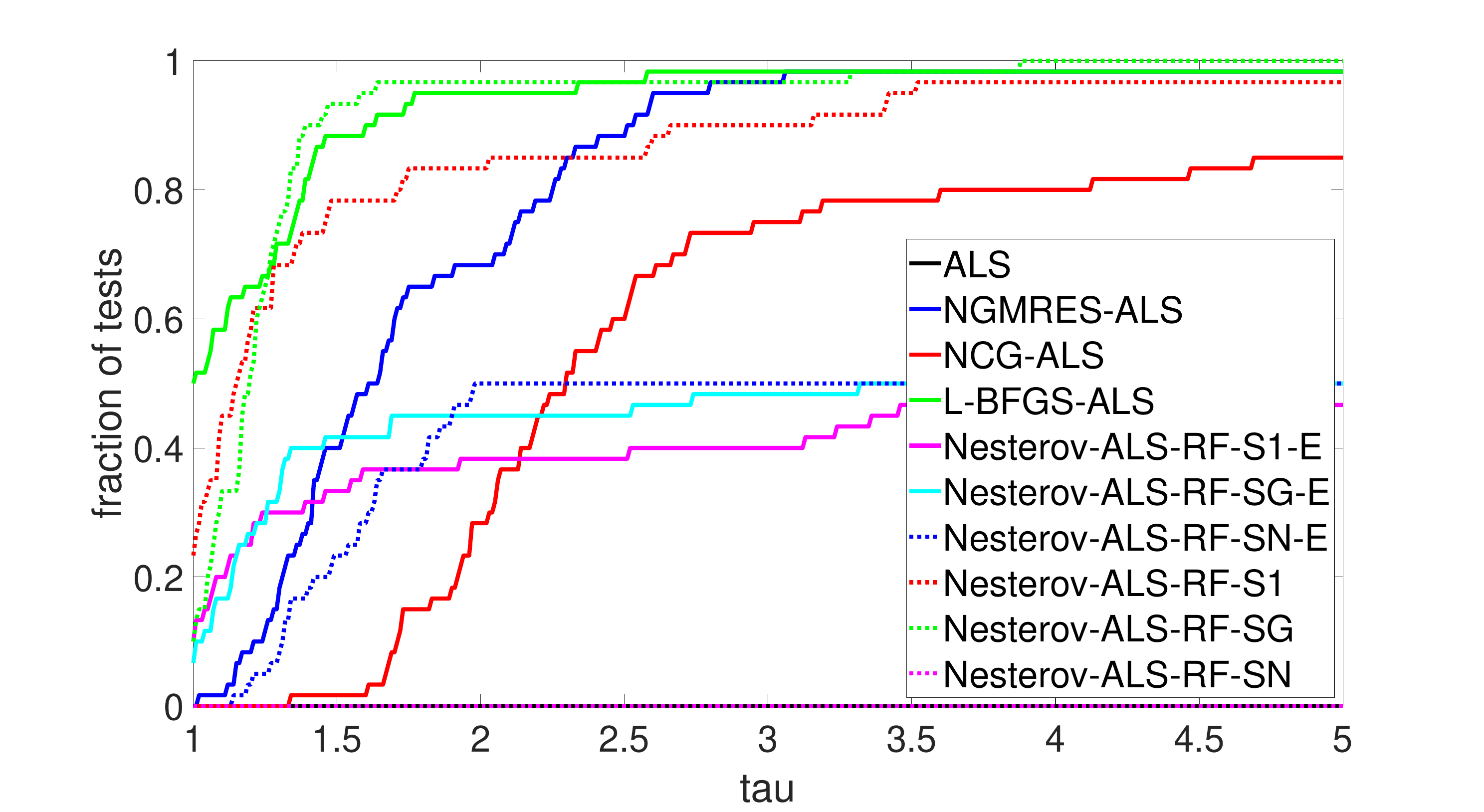}}
		\caption{variants without delay}
	\end{subfigure}
	\caption{Synthetic test problems. $\tau$-plots comparing variants of function restart.}
	\label{fig:fr}
	\end{centering}
\end{figure*}

\begin{figure*}[h!]
	\begin{centering}
	\begin{subfigure}{.8\linewidth}
		\centering{\includegraphics[width=.87\linewidth]{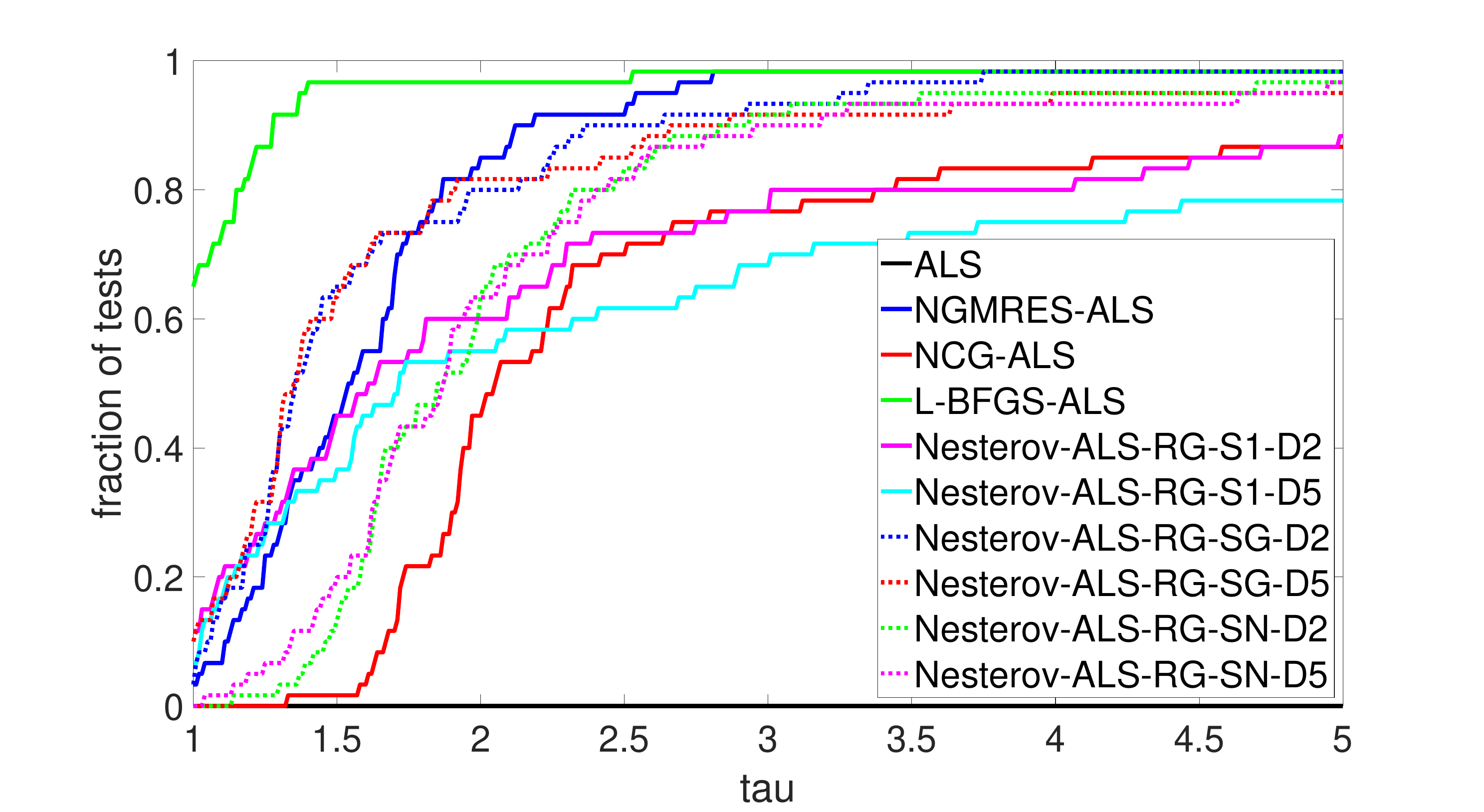}}
		\caption{variants with delay}
	\end{subfigure}

	\begin{subfigure}{.8\linewidth}
		\centering{\includegraphics[width=.87\linewidth]{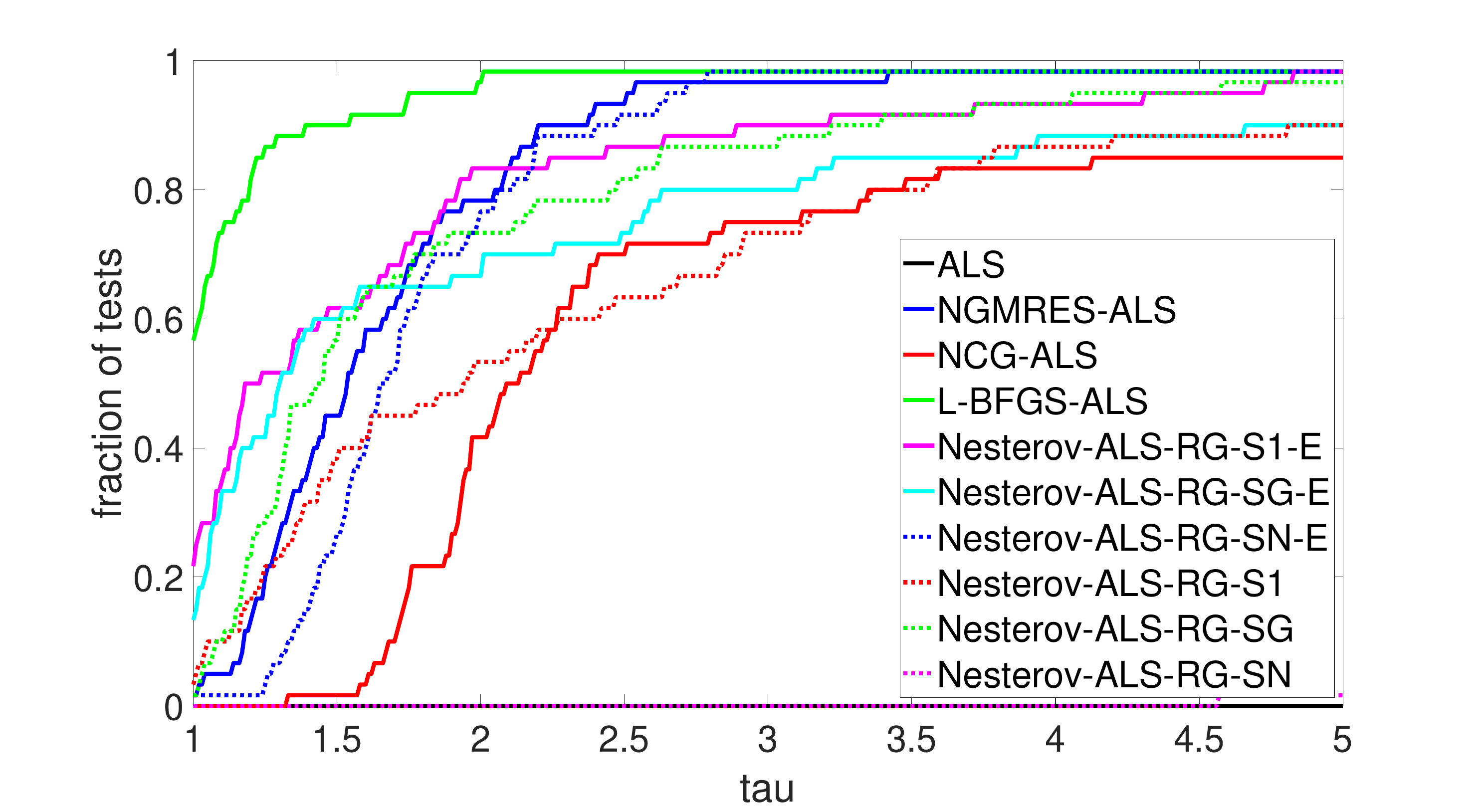}}
		\caption{variants without delay}
	\end{subfigure}
	\caption{Synthetic test problems. $\tau$-plots comparing variants of gradient restart.}
	\label{fig:gr}
	\end{centering}
\end{figure*}

\begin{figure*}[h!]
	\begin{centering}
	\centering{\includegraphics[width=\linewidth]{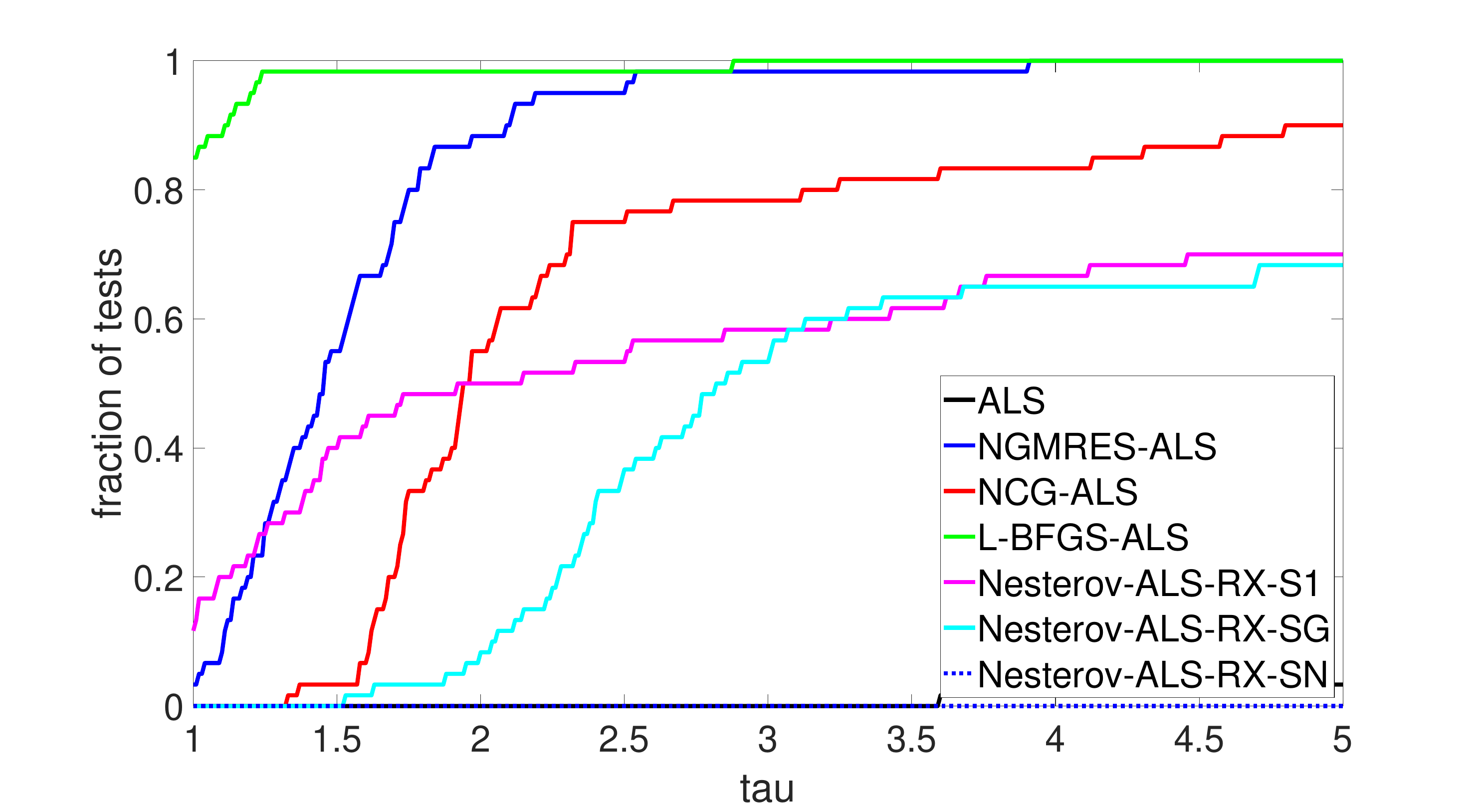}}
	\caption{Synthetic test problems. $\tau$-plot comparing variants of speed restart.}
	\label{fig:sr}
	\end{centering}
\end{figure*}

\Cref{fig:fr,fig:gr,fig:sr} show $\tau$-plots for variants of the restarted Nesterov-ALS schemes, for the case of function restart (RF, \Cref{fig:fr}), gradient restart (RG, \Cref{fig:gr}), and speed restart (RX, \Cref{fig:sr}), applied to the synthetic test problems.

For each of the restart mechanisms, several of the restarted Nesterov-ALS variants typically outperform ALS, NCG-ALS\cite{sterck2015nonlinearly} and NGMRES-ALS\cite{sterck2012nonlinear}.

Several of the best-performing restarted Nesterov-ALS variants are also competitive with the best existing
nonlinear acceleration method for ALS we compare with, LBFGS-ALS\cite{sterck2018nonlinearly}, and they are much easier to implement. 

Among the restart mechanisms tested, function restart (\Cref{fig:fr}) substantially outperforms gradient restart
(\Cref{fig:gr}), and, in particular, speed restart (\Cref{fig:sr}).

The $\tau$-plots confirm that Nesterov-ALS-RF-SG, using function restarting and gradient ratio momentum
weight, consistently performs as one of the best methods, making it our recommended choice for ALS acceleration.

%%%%%%%%%%%%%%%%%%%%%%%%%%%%%%%%
%%%%%%%%%%%%%%%%%%%%%%%%%%%%%%%%
%%%%%%%%%%%%%%%%%%%%%%%%%%%%%%%%
\end{document}